\documentstyle[12pt]{article}

\newtheorem{lemma}{Lemma}
\newtheorem{theorem}[lemma]{Theorem}

\newtheorem{corollary}[lemma]{Corollary}
\newtheorem{proposition}[lemma]{Proposition}

\newcommand{\eop}{\hfill $/$\hspace*{-.1cm}$/$\hspace*{-.1cm}$/$\vspace{.1in}}

\newcommand{\subnum}[1]{
\begin{center}
{\bf #1}
\end{center}
}

\newcommand{\Ff}{{\cal F}}
\newcommand{\Gg}{{\cal G}}

\begin{document}

\section*{A Giraud-type characterization of the simplicial categories associated
to closed model categories as $\infty$-pretopoi}

\noindent
Carlos Simpson
\newline
CNRS, UMR 5580, 31062 Toulouse CEDEX

\bigskip

In SGA 4 \cite{SGA4}, one of the principal building blocks of the theory of
{\em topoi} is Giraud's theorem, which says that the condition of a
$1$-category $A$ being the category of sheaves on a site, may be
characterized by
intrinsic, internal conditions in $A$. The intrinsic conditions are basically
existence of certain limits and colimits, plus a condition about
generation by a small set of objects.

In this paper, we will present a generalization of this theorem to the
situation of {\em simplicial categories} (by which we mean simplicially
enriched categories) or equivalently {\em Segal categories} (\cite{DKS}
\cite{SimpsonFlexible}, \cite{SimpsonEffective}), or {\em complete Segal spaces}
(Rezk \cite{Rezk}). One can easily imagine generalizing the internal conditions
of existence of limits or colimits (these become conditions of existence of
homotopy limits or colimits). On the other hand, the condition which we take as
a generalization of the condition of being the category of sheaves on a site, is
the condition of coming from a {\em closed model category} \cite{Quillen}.
Recall that Dwyer and Kan associate to any closed model category $M$ its {\em
simplicial localization} $L(M)$ which is a simplicial category
\cite{DwyerKan}. If $M$ is a {\em
simplicial closed model category} in the sense of Quillen, then $L(M)$ is
equivalent to the simplicial category of fibrant and cofibrant objects of $M$.
It is this simplicial category $L(M)$ which represents the {\em homotopy theory}
(including information about all higher-order homotopies) which comes out of
$M$.

We attack the very natural
question of characterizing which simplicial categories $A$ are equivalent to
ones of the form $L(M)$ for closed model categories $M$.
This formulation of the question is closely related to some of the entries in
the ``Model Category'' section of
M. Hovey's recent ``problem list''
\cite{HoveyList}.

The first, easy but fundamental observation is that if $M$ is a closed model
category (admitting all small limits and colimits as it is now customary to
assume), then $L(M)$ admits small homotopy limits and colimits. In particular,
not every simplicial category will be equivalent to one of the form $L(M)$.
Our characterization is that this necessary condition is basically sufficient;
however, one has to add in an additional set-theoretic hypothesis about
small generation which in practice will always hold.  This first easy
observation
came from thinking about C. Rezk's terminology of calling his version of the
closed model category of Segal categories, the ``homotopy theory of homotopy
theories''.

Our answer is, as stated above, analogous to Giraud's theorem. To be quite
precise, the analogy is not complete. In effect, the internal conditions on $A$
which come out are existence of homotopy colimits, and small generation. These
turn out to imply existence of homotopy limits; however one does not get any
sort of exactness properties allowing one to commute limits and colimits, and
indeed one can find examples of closed model categories $M$ such that
$L(M)$ does not have these exactness properties. Thus, in the statement of our
theorem, we refer to our equivalent conditions as defining a notion of {\em
$\infty$-pretopos}, and reserve the name {\em $\infty$-topos} for an
$\infty$-pretopos satisfying additional exactness properties.

Another remark is that we are not able to treat {\em all} closed model
categories, nor does this seem natural in the context of a Giraud-type theorem.
Rather we speak only of {\em cofibrantly generated closed model categories}
see \cite{DHK} \cite{Hovey} \cite{Hirschhorn}. Almost all known closed model
categories (here as usual we only consider ones in
which all small limits and colimits exist) are cofibrantly generated.

Hovey also states in
\cite{HoveyList} that D. Dugger has shown that any cofibrantly generated closed
model category is Quillen-equivalent to a simplicial one; thus the reader of the
present introduction who is unfamiliar with Dwyer-Kan may assume that we are
speaking of simplicial model categories and may replace $L(M)$ by the simplicial
category $M_{cf}$ of cofibrant and fibrant objects.

Here is a shortened version of the statement. As a matter of notation, we speak
in the introduction of ``simplicial categories''; the notion of {\em
equivalence} is that which was explored by Dwyer and Kan \cite{DwyerKan}.
This is just the obvious notion of ``fully faithful and essentially surjective''
where ``fully faithful'' means inducing weak equivalences of simplicial $Hom$
sets, and ``essentially surjective'' means essential surjectivity of the
truncated morphism on homotopy $1$-categories. However, with this definition
an {\em equivalence} between two simplicial categories means a string of
functors which are equivalences, possibly going in different directions.
See below for a bit more explanation.  We also refer to the body of the paper
for the definitions of homotopy colimit, generation by homotopy colimits, and
smallness.

\begin{theorem}
\label{one}
(cf Theorem \ref{main} p. \pageref{main})
Suppose $A$ is a simplicial category.  The following conditions are equivalent:
\newline
(i)\,\,  There is a cofibrantly  generated closed model category $M$ such that
$A$ is equivalent to the Dwyer-Kan simplicial localization $L(M)$;
\newline
(ii) \,\, $A$ admits all small homotopy colimits, and there is a small subset of
objects of $A$ which are $A$-small, and which generate $A$ by homotopy
colimits.
\end{theorem}

We call a simplicial category satisfying the conditions of the theorem, an {\em
$\infty$-pretopos}. If in addition a certain exactness condition is satisfied
(see the statement of Theorem \ref{main} for details) then we say that $A$ is
an {\em $\infty$-topos}.

The possibility of having a reasonable notion of
$n$-topos was predicted in \cite{SimpsonRealization}. This
prediction came about due to the influence of correspondence with C. Teleman who
at the time was telling me about pullbacks of simplicial presheaves under
morphisms of sites.  Of course, like most of what we do here, this idea is very
present in spirit throughout \cite{Grothendieck}.

A word about rigour and level of detail in this version of the present
paper. At several places in the argument, we skip verification of some details.
These are mostly details concerning ``homotopy-coherent category theory''
as done with Segal categories. They are all generalizations to the
``weak-enriched'' setting of classical statements in category theory, so it
seems completely clear that the statements in question are true. It also seems
clear that in the relatively near future, techniques will have sufficiently
advanced in order to cover these questions. Finally, it seems likely that using
some of the other approaches (such as Cordier-Porter \cite{CordierPorter} or
the model category of Dwyer-Hirschhorn-Kan \cite{DHK}), a significant number of
these details could be verified relatively easily---the reason I haven't taken
that route is lack of familiarity with those approaches. However, at the time of
writing of the present version, I have not verified the details any further than
what is written down below. One could say that the present paper is premature
in this sense, but the result seemed interesting enough to justify writing it
up quickly. In order to clarify  matters, the places where this problem occurs
are marked with the symbol $(\otimes )$.

{\em Acknowledgements:} I would like to thank very much M. Hovey, C. Rezk, A.
Hirschowitz, P. Hirschhorn, C. Teleman, Z. Tamsamani, B. Toen, and J. Tapia, for
their important contributions to the realization of the idea  outlined in the
present paper.

\subnum{Segal categories}

A {\em simplicial category} means a category $C$ enriched over simplicial sets.
In other words, for every pair of objects $x,y\in ob(C)$, we have a simplicial
set $Hom _C(x,y)$. In order to conform with our notations for Segal categories,
we shall denote this simplicial set by $C_{1/}(x,y)$. In the case of  a
simplicial category, composition of morphisms is a map of simplicial sets
$$
C_{2/}(x,y,z):= C_{1/}(x,y) \times C_{1/}(y,z) \rightarrow C_{1/}(x,z),
$$
and this operation is strictly associative. In view of the strict associativity,
we obtain a bisimplicial set (i.e. simplicial simplicial set) by setting
$$
C_{p/}:= \coprod _{x_0,\ldots , x_p} C_{p/}(x_0,\ldots , x_p),
$$
with
$$
C_{p/}(x_0,\ldots , x_p) := C_{1/}(x_0,x_1) \times \ldots \times
C_{1/}(x_{p-1}, x_p).
$$
Here we set $C_0:= C_{0/}:= ob(C)$. This bisimplicial set has the property that
for any $m$, the ``Segal map''
$$
C_{m/} \rightarrow C_{1/} \times _{C_0} \ldots \times _{C_0}C_{1/}
$$
is an isomorphism; and conversely any bisimplicial set such that the
simplicial set $C_{0/}$ is a discrete set  which we denote $C_0$ or $ob(C)$, and
such that the above Segal maps are isomorphisms, corresponds to a simplicial
category.
The composition is obtained by using the third face map $C_{2/} \rightarrow
C_{1/}$.

The above presentation of the notion of ``simplicial category'' motivates the
definition of ``Segal category''---a {\em Segal category} is just a
bisimplicial set such that the simplicial set $C_{0/}$ is a discrete set which
we denote $C_0$ or $ob(C)$, and such that the Segal maps are weak equivalences
of simplicial sets. A simplicial category thus gives rise to a
Segal category, and we shall sometimes call the Segal categories which arise
in this way ``strict''.

Suppose $C$ is a Segal category. As suggested by the previous notation, for any
sequence of objects $x_0,\ldots , x_p\in ob(C)$ we obtain a simplicial set
$C_{p/}(x_0,\ldots , x_p)$ defined as the inverse image of $(x_0,\ldots , x_p)$
under the map (given by the $p+1$ ``vertex'' maps)
$$
C_{p/} \rightarrow C_0 \times \ldots \times C_0.
$$
We think of the simplicial set $C_{1/}(x,y)$ as being the space of maps
from $x$ to $y$ in $C$.
The Segal condition can be rewritten as saying that the morphism
(given by the $p$ ``principal edge'' maps)
$$
C_{p/}(x_0,\ldots , x_p)\rightarrow
C_{1/}(x_0,x_1) \times \ldots \times C_{1/}(x_{p-1}, x_p)
$$
is a weak equivalence. In particular, the ``composition of morphisms in
$C$'' is
given by the diagram
$$
C_{1/}(x,y)\times C_{1/}(y,z) \stackrel{\cong}{\leftarrow}
C_{2/}(x,y,z)\rightarrow C_{1/}(x,z).
$$

The notion of Segal category is based in an obvious way on Segal's weakened
notion of ``topological monoid'' \cite{Segal} (which is the case where $ob(C)$
contains only one element), although Segal himself never seems to have written
anything suggesting to look at this notion for several objects. This notion
{\em per se} first appears in Dwyer-Kan-Smith \cite{DKS} where they also show
the equivalence between Segal categories-up-to-equivalence and simplicial
categories-up-to-equivalence (see below).

This notion later appeared in an
{\em ad hoc} way in my preprint \cite{SimpsonFlexible} (I was
unaware of \cite{DKS} at the time and until fairly recently); and it appears as
the basic idea which is iterated in Tamsamani's definition of weak $n$-category
\cite{Tamsamani}.  Further occurences are in my preprint \cite{SimpsonEffective}
and the joint paper \cite{HirschowitzSimpson}.

A couple of closely
related notions are use by Rezk in \cite{Rezk}. He defines a
notion of {\em Segal space} which is a simplicial set satisfying the
condition that the Segal maps are equivalences but not necessarily the
condition that $C_{0/}$ be a discrete set; consequently he includes a ``Reedy
fibrant'' condition in the definition in order to make sure that the fiber
products involved in the definition of the Segal maps are homotopically
correct ones. He also defines a notion of {\em complete Segal space}
which basically says that the simplicial set $C_{0/}$ should itself
correspond to the space which is the realization of the subcategory obtained
by only looking at invertible (up-to-homotopy) morphisms in $C$. We will
state without proof below what should be the relation between Rezk's notions and
our own.

There are also other related notions such as various notions of {\em
$A_{\infty}$-category} see for example Batanin \cite{Batanin};
and more generally there are several definitions of $n$-category alternative
to Tamsamani's definition and which should also have variants for weak
simplicial categories, see Baez-Dolan \cite{BaezDolan} for example. These
other notions should be directly related to our own but we don't go into that
here.

Finally, we note that the above notions should be viewed as substitutes
for the notion of ``$1$-groupic $\infty$-category'' i.e. an $\infty$-category
in which the $i$-morphisms are invertible (up to equivalence) for $i\geq 2$.
We leave it to the reader to make this notational translation.

We shall use the framework of ``Segal categories'' throughout the rest of the
paper, although we sometimes speak of the relationship with the
classical notion of simplicial  category.  The reader is refered to
\cite{SimpsonLimits} and \cite{HirschowitzSimpson} for any further details
and introductory materiel that we may leave out in our brief discussion which
follows.

By abuse of notation, we may sometimes forget to put in the qualifier ``Segal''
and just use the word ``category'' for ``Segal category''. In order to avoid
confusion, we will try to systematically use the terminology {\em $1$-category}
for classical (non-simplicial) categories.

A morphism of Segal categories
$C\rightarrow D$ is said to be {\em fully faithful} if for every $x,y\in
ob(C)$, the morphism $C_{1/}(x,y)\rightarrow  D_{1/}(x,y)$ is a weak
equivalence of simplicial sets. This is the natural generalization of the
corresponding notion in category theory; however one should be careful that the
separate notions of ``full'' and ``faithful'' don't have reasonable
generalizations to the present theory, because there is no way of decomposing
the condition of being a weak equivalence of simplicial sets, into
``injectivity plus surjectivity''.  For this reason, huge swaths of the
argumentation which is employed in SGA 4 \cite{SGA4} are no longer available and
we are forced to look for more intrinsic reasoning.

If $C$ is a Segal category, define a $1$-category denoted $ho(C)$ with the same
objects as $C$, by setting
$$
ho(C)_{1/}(x,y):= \pi _0(C_{1/}(x,y)).
$$
We say that a morphism of Segal categories $C\rightarrow D$ is {\em essentially
surjective} if the resulting morphism of $1$-categories
$$
ho(C)\rightarrow ho(D)
$$
is essentially surjective. We say that a morphism of Segal categories is an
{\em equivalence} if it is fully faithful and essentially surjective.

In the context of simplicial categories this notion of equivalence was
introduced by Dwyer and Kan in \cite{DwyerKan}. (A morphism of simplicial
categories is an equivalence if and only if the corresponding morphism of Segal
categories is an equivalence.) In the context of $n$-categories this notion was
called ``external equivalence'' by Tamsamani in \cite{Tamsamani}. In his
situation of Segal spaces, this notion was called ``Dwyer-Kan equivalence'' by
Rezk in \cite{Rezk}.

We say
that a morphism in $C$ (i.e. a vertex of $C_{1/}(x,y)$) is an {\em equivalence}
if its image in $ho(C)_{1/}(x,y)$  is an isomorphism in $ho(C)$. This
corresponds to what Tamsamani called ``internal equivalence'' in
\cite{Tamsamani}. The essential surjectivity condition can be expressed as
saying that every object of $D$ is equivalent (in this ``internal'' sense) to
an object coming from $C$.

We often use the terminology {\em full subcategory} for a fully faithful functor
of Segal categories $C\rightarrow D$ which is injective on objects.  In
this case, up to equivalence in the variable $C$, we may assume that the
morphism is actually an isomorphism on all of the  $C_{p/} (x_0,\ldots , x_p)$.
With this convention, the intersection of full subcategories is again a
full subcategory. Furthermore, we say that a full subcategory $C\subset D$
is {\em saturated} if it satisfies the ``saturation condition'' that whenever
$x\in ob(C)$ and $y$ is (internally) equivalent to $x$, then $y\in ob(C)$ too.
Again, the intersection of saturated full subcategories is again a saturated
full subcategory.

\subnum{Strictification}

We can now explain the comparison result of Dwyer, Kan, Smith which was alluded
to above. Let ${\rm splCat}$ denote the $1$-category of simplicial categories,
and let ${\rm SegCat}$ denote the $1$-category of Segal categories. Let
$$
Ho({\rm splCat}) (\mbox{resp.} \;\;\; Ho({\rm splCat}))
$$
denote the Gabriel-Zisman
localizations of these categories by inverting the equivalences.
We say that two simplicial categories (or two Segal categories) are {\em
equivalent} if they project to isomorphic objects in these homotopy categories.
Dwyer, Kan and Smith in the last few pages of \cite{DKS} show that the
morphism
$$
Ho({\rm splCat})\rightarrow Ho({\rm SegCat})
$$
is an equivalence  of categories. Among other things, this says that any Segal
category can be ``strictified'', i.e. made equivalent (in the above sense) to a
simplicial category. We should take this occasion to stress that, as $Ho({\rm
splCat})$ and $Ho({\rm splCat})$ are Gabriel-Zisman localizations, one can have
two objects (simplicial categories or Segal categories) which are equivalent
but without there being any
actual morphism between the two; the ``equivalence'' in question might be
realizable only as a chain of morphisms which are equivalences, going in
different directions. This situation is improved by the introduction of closed
model structures as we shall explain below (and in particular if ever it is
necessary to go through a chain of equivalences, at least one can restrict to
looking at chains of length $2$). In the statement of Theorem \ref{one}, it is
the present notion of equivalence which is used.

In view of the strictification result of Dwyer-Kan-Smith, we may at many places
in the present paper assume that the Segal categories we are dealing with are
actually simplicial categories. This can simplify the problem of composing
morphisms and the like.

\subnum{Closed model structures}

There are several possible closed model structures which can be used to
attack the homotopy category $Ho({\rm
splCat})\cong Ho({\rm splCat})$.
What seems to be historically the first
is that of Dwyer-Hirschhorn-Kan \cite{DHK}.
\footnote{The draft of \cite{DHK} that I have is dated at approximately the same
time as \cite{SimpsonVK} but earlier versions of \cite{DHK} had apparently
been in limited circulation for some time.}

To introduce the structure of \cite{DHK}, we first
point out that Dwyer and Kan obtained (essentially trivially) a closed model
structure  on ${\rm splCat}$ in \cite{DwyerKan} where the weak equivalences were
the equivalences which induce isomorphisms on objects. In this structure,
the fibrations are the
morphisms of simplicial categories which induce fibrations of the individual
simplicial $Hom$ sets. The cofibrations are closely related to the free
resolutions which are used throughout \cite{DwyerKan}, and the cofibrant
objects are just the simplicial categories which are free at each stage. This
closed model structure is not the one which we are actually interested in
(although it can be useful in a preliminary way), because we are interested in
understanding the equivalences which are essentially surjective but not
isomorphisms on objects. This problem was rectified in \cite{DHK} where a closed
model structure on  ${\rm splCat}$ is given, with the following properties.
The cofibrations are the same as in the previous structure; and the weak
equivalences are the ``Dwyer-Kan equivalences'' as described above. This leads
to a more restrictive notion of fibration than that which occurs in their first
structure. However, the fibrant objects are the same as in the previous
structure, namely the simplicial categories $C$ with $C_{1/}(x,y)$ being
fibrant simplicial sets. To sum up what is going on  here, we can say that in
order to correctly calculate the morphisms between two simplicial categories
$C$ and $D$, one must make a replacement $D\rightarrow D^f$ by an equivalent
one in which the simplicial $Hom$-sets are fibrant, and one must make a
replacement $C^c\rightarrow C$ with $C^c$ cofibrant (which essentially means
taking a free resolution). Now $Hom(C^c, D^f)$ contains representatives for all
of the homotopy classes of morphisms from $C$ to $D$ in $Ho({\rm splCat})$.

The main drawback of the Dwyer-Hirschhorn-Kan closed model structure is that
the cofibrant replacement $C^c\rightarrow C$ is not compatible with direct
product. Thus one does not obtain (in any direct way) an internal
$\underline{Hom}(C,D)$. This internal $\underline{Hom}$ will be crucial for the
arguments in the present paper.

In \cite{SimpsonVK} is given a closed model structure for $n$-categories.
This is essentially the same problem as for Segal categories, and in
\cite{SimpsonEffective}, the closed model structure for Segal
categories was announced with the statement that the proof is the same as in
\cite{SimpsonVK}.  A complete proof was written up in
\cite{HirschowitzSimpson}. This closed model structure yields as underlying
homotopy category $Ho({\rm SegCat})$, and it is ``internal'' i.e. admits a
homotopically correct internal $\underline{Hom}$. We shall use this structure in
the present paper.

Before getting to a more detailed description, we note that Rezk constructs a
closed model structure for what he calls {\em complete Segal spaces} in
\cite{Rezk}. Rezk's closed model structure again yields
as underlying homotopy category a category which is equivalent to
$Ho({\rm SegCat})\cong Ho({\rm SegCat})$ (this fact follows immediately from
the statements in \cite{Rezk} plus the strictification result
of Dwyer-Kan-Smith). And Rezk's closed model structure is ``internal'',
in other words it can be used to calculate $\underline{Hom}(C,D)$.
Thus it should be possible to write the present paper using Rezk's structure
rather than my own.  The obvious conjecture is that Rezk's structure and my own
are Quillen-equivalent. We don't prove this here, but it is probably an easy
consequence of everything that is said in Rezk's preprints and my own---the
only problem being to digest all of that!

We should also at this point mention another approach, which is the
``homotopy-coherent'' approach of Cordier and Porter \cite{CordierPorter}.
They define a simplicial category $Coh(C,D)$ for any two simplicial categories
$C$ and $D$. This should be equivalent to the $\underline{Hom}$ constructed
in either my model category or Rezk's model category, and again it should be
possible to write the present paper using Cordier-Porter's theory (and indeed
this might be advantageous in many places).

Invoking the principle that the author of a paper is allowed to choose which
approach he wants to use, we will use the closed model structure of
\cite{SimpsonVK},
\cite{SimpsonEffective} and
\cite{HirschowitzSimpson}, which we now describe. The first step is to define
a $1$-category of {\em Segal precats} denoted $SePC$. The objects are the
bisimplicial sets (denoted as above $p\mapsto X_{p/}$ with $X_{p/}$ denoting
a simplicial set) such that $X_{0/}$ is a discrete set which we denote $C_0$ or
$ob(X)$. The morphisms in $SePC$ are just morphisms of bisimplicial sets.
This category admits all small limits and colimits. We define the {\em
cofibrations} to be the monomorphisms in this category, in other words the
injections of bisimplicial sets. It remains to be seen how  to define the
weak equivalences. For this, note that the category ${\rm SegCat}$ is a
subcategory of $SePC$. The main step (we refer to \cite{SimpsonVK},
\cite{SimpsonEffective} and
\cite{HirschowitzSimpson} for the details of which) is an essentially unique
``projection functor''
$$
SeCat: SePC \rightarrow {\rm SegCat}\subset SePC,
$$
together with a natural transformation $\eta _X: X\rightarrow SeCat(X)$,
such that $\eta _X$ is an equivalence if $X$ is already a Segal category.
This is a variant of the well-known notion of ``monad'' in category theory,
a variant which uses the notion of equivalence (rather than isomorphism) in the
target subcategory ${\rm SegCat}$. We think of $SeCat(X)$ as being the Segal
category generated by the ``generators and relations'' $X$. In
\cite{SimpsonEffective} the operation $X\mapsto SeCat(X)$ is analyzed explicitly
and shown to have good effectivity properties.

Now we say that a morphism $X\rightarrow Y$ is a {\em weak equivalence} if the
resulting morphism of Segal categories $SeCat(X)\rightarrow SeCat(Y)$ is
an equivalence in the sense explained above. This gives rise to the notion of
{\em trivial cofibration} (a cofibration which is a weak equivalence) and hence
to the notion of {\em fibration} (a morphism which satisfies lifting for all
trivial cofibrations). It is shown in \cite{SimpsonVK}
and  \cite{HirschowitzSimpson} that $SePC$ with these three classes of
morphisms is a cofibrantly generated closed model category. One thing to note
is that the fibrant objects of $SePC$ are themselves Segal categories, i.e.
$$
SePC _f \subset {\rm SegCat}.
$$
It follows that
$$
Ho (SePC ) \cong Ho (SePC_f) \cong Ho({\rm SegCat}).
$$

The closed model category $SePC$ is ``internal'', see
\cite{SimpsonVK} and  \cite{HirschowitzSimpson}. This basically means that
the cartesian product is a monoidal structure in the sense of Hovey {\em et
al.}. The effect of this property is that we have a notion of internal
$\underline{Hom}$ in $SePC$. This is defined by the adjunction property that
for any Segal precat $E$, a morphism
$$
E\rightarrow \underline{Hom}(A,B)
$$
is the same thing as a morphism (in $SePC$)
$$
A\times E \rightarrow B.
$$
Now if $B$ is a fibrant Segal category (i.e. a fibrant object in $SePC$) then
for $A$ any Segal precat, $\underline{Hom}(A,B)$ is again a fibrant Segal
category. In the case where the second variable is fibrant, formation of the
internal $\underline{Hom}$ is compatible with weak equivalences in both
variables. We will make heavy use of this internal $\underline{Hom}$, bearing
in mind that whenever it is used, the second variable has to be made fibrant.

The above discussion leads to the notion of {\em natural transformation}
between two functors of Segal categories. If $A$ and $B$ are Segal categories
(with $B$ assumed to be fibrant) and if $f,g: A\rightarrow B$ are morphisms,
a {\em natural transformation from $f$ to $g$} is a vertex of the simplicial set
$$
\eta \in \underline{Hom}(A,B)_{1/}(f,g).
$$
In general for a Segal category $C$, a vertex of $C_{1/}(x,y)$ is the same
thing as a morphism $I\rightarrow C$ (where $I$ is the $1$-category with two
objects $0,1$ and one arrow $0\rightarrow 1$) such that $0$ goes to
$x$ and $1$ goes to $y$. Apply this with $C=\underline{Hom}(A,B)$. We get that
a natural transformation from $f$ to $g$ is the same thing as a morphism
$$
\eta : A\times I \rightarrow B
$$
such that $\eta |_{A\times 0}=f$ and $\eta |_{A\times 1} = g$.

The internal $\underline{Hom}$ is used in \cite{HirschowitzSimpson} (following
the same idea in the case of $n$-categories in \cite{SimpsonVK}) to define the
{\em Segal $2$-category $1SeCAT$ of all Segal categories}. This has for objects
the fibrant Segal categories, and between two objects $A,B$ one takes as Segal
category of morphisms the internal $\underline{Hom}(A,B)$. We get a strict
category enriched over fibrant Segal categories, which yields a Segal
$2$-category. We refer to \cite{SimpsonVK} and  \cite{HirschowitzSimpson} for
more details; this will not be used in the remainder of the present paper.

We now indicate a sketch of how one should obtain the relationship between
the above closed model category and Rezk's closed model category \cite{Rezk}
of complete Segal spaces which we shall denote $RC$ for the
present discussion.  If $A$ is  a Segal category, let $rf(A)$ be a Reedy-fibrant
replacement of $A$ as bisimplicial set. Then $rf(A)$ is a Segal space in Rezk's
terminology. Now Rezk has a construction which replaces a Segal space by a
complete Segal space, which we will denote by $crf(A)$. This gives a functor
going from the category of Segal categories to the category of complete Segal
spaces. It descends to the Gabriel-Zisman (or even Dwyer-Kan) localizations
where we divide out by equivalences (Rezk states that his construction takes
Dwyer-Kan equivalences of Segal spaces, to equivalences of complete Segal
spaces). In the other direction, given a complete Segal space $X$, we can
discretize the space of objects and chop up the other spaces accordingly (in the
minimal way so that the transition morphisms remain continuous). This yields a
Segal category. Again, this construction takes equivalences to equivalences.
Thus we obtain an equivalence of $1$-categories between the homotopy category of
Segal categories, and the homotopy category of complete Segal spaces:
$$
Ho({\rm SegCat} )\cong Ho(SePC) \cong Ho(RC).
$$
Furthermore, on the level of Dwyer-Kan localizations we obtain an equivalence of
simplicial categories
$$
L(SePC) \cong L(RC).
$$
Technically speaking, there is probably some remaining
verification to be done here, for example verifying that the two constructions
are really inverses.
It would also be nice to set up a Quillen  equivalence
between the two  model categories, and to verify that the equivalences are
compatible with internal $\underline{Hom}$.

This last compatibility is already
obtained on a homotopy-theoretic level in the following way: it was observed
(e.g. in \cite{HirschowitzSimpson}) that  if a closed model category $M$ is
``internal'', then its Dwyer-Kan localization $L(M)$ is a simplicial category
admitting internal $Hom$ as defined in an appropriate way. In this case, the
internal $\underline{Hom}(X,Y)$ (for $X,Y\in L(M)$) may be characterized in a
way which is internal to $L(M)$. This applies both to $SePC$ and to Rezk's
closed model category $RC$. Since the two localizations $L(SePC)$ and
$L(RC)$ are
equivalent (by the argument sketched above), this shows that the internal
$\underline{Hom}(X,Y)$ are equivalent in $L(SePC)$ and $L(RC)$.  Another way to
recast this remark is to point out that, applying the result of Dwyer-Kan-Smith
\cite{DKS} we obtain an equivalence with the Dwyer-Kan localization of the
Dwyer-Hirschhorn-Kan model category (we denote the latter by $DHK$)
$$
L(DHK)\cong  L({\rm splCat})\cong L({\rm SegCat})\cong L(SePC) \cong L(RC),
$$
and existence of the internal closed model categories $SePC$ and Rezk's $RC$
can be viewed as ways of proving that the simplicial category $L(DHK)$
admits an  internal $\underline{Hom}$.

We close this subsection on a slightly more technical note.
In many places, the notation $\Upsilon$  introduced in
\cite{SimpsonLimits} is crucial for correctly manipulating Segal categories in
our method. We refer to there (or to any of a number of my more recent
preprints where this notation is used) for details and examples. A rapid
overview
would say that if $E$ is a simplicial set then we obtain a Segal precat
$\Upsilon (E)$  having two objects denoted $0,1$, and having $E$ as simplicial
set of morphisms from $0$ to $1$. In the case $E=\ast$ we recover $\Upsilon
(\ast )=I$, the $1$-category with objects $0$ and $1$ and a single morphism
$0\rightarrow 1$. This has a sort of universal property: for any Segal precat
$A$, a morphism $E\rightarrow A_{1/}(x,y)$ is the same thing as a morphism
$$
\Upsilon (E)\rightarrow A
$$
sending $0$ to $x$ and $1$ to $y$.

More generally if $E,F$ are simplicial sets then we obtain $\Upsilon ^2(E,F)$
which has objects $0,1,2$ and $E$ as morphisms from $0$ to $1$; $F$ as
morphisms from $1$ to $2$; and $E\times F$ as morphisms from $0$ to $2$.
This latter is useful for dividing up a square into two triangles: one has the
pushout formula
$$
\Upsilon (E)\times \Upsilon (F) \cong \Upsilon ^2(E,F) \cup ^{\Upsilon (E\times
F)} \Upsilon ^2(F,E).
$$
Finally, the existence of weak compositions is manifested in the statement that
the inclusion
$$
\Upsilon (E) \cup ^{\{ 1\} } \Upsilon (F) \rightarrow \Upsilon ^2(E,F)
$$
is a trivial cofibration.

\subnum{Simplicial sets and cartesian families}

Let $S$ denote the simplicial category of all fibrant simplicial sets.
It has for objects the fibrant simplicial sets $K$, and for simplicial $Hom$
sets the internal $\underline{Hom}(K,L)$ of simplicial sets.

Unfortunately, $S$ is not fibrant as a Segal category. Thus we must fix a
fibrant replacement $S\rightarrow S'$ (i.e. an equivalence of Segal categories
with $S'$ fibrant). Note here that $S'$ cannot be a strict simplicial category.
This fibrant replacement is a source of most of the technical difficulties
which were encountered in \cite{SimpsonLimits} and
\cite{HirschowitzSimpson}. The best way to get around these problems, at least
in the context of the theory we are exposing here, is the canonical fibrant
replacement defined using the notion of ``cartesian family'' in
\cite{SimpsonAspects}. This was constructed in the context of
$n$-categories, giving a fibrant replacement for the $n+1$-category $nCAT$ of
all $n$-categories. We describe here the variant for obtaining a fibrant
replacement for $S$ (note that in the notation of \cite{HirschowitzSimpson}, a
simplicial set is a Segal $0$-category and $S=0SeCAT$; the variant we are about
to describe is obtained from the discussion in \cite{SimpsonAspects}
by  substituting ``$0Se$'' for ``$n$'').

For ease of use in the rest of the paper, we consider ``contravariant''
cartesian families; these will correspond to functors $A^o\rightarrow S'$,
and this constitutes a change with respect to \cite{SimpsonAspects} where
``covariant'' cartesian families were considered.

Suppose $A$ is a Segal category, considered as a bisimplicial set. A {\em
(contravariant)
precartesian family (of simplicial sets) over $A$} is a morphism of bisimplicial
sets
$$
\Ff \rightarrow A
$$
satisfying the  ``cartesian property'' which we now explain. We first establish
some notations: $\Ff _{p/}$ is the simplicial set obtained by
putting $p$ in the first bisimplicial variable; thus $\Ff _{p/} \rightarrow
A_{p/}$. For objects $x_0, \ldots , x_p\in ob(A)$, we denote by
$$
\Ff _{p/}(x_0,\ldots , x_p)
$$
the inverse image of $A_{p/}(x_0,\ldots , x_p)$. It is also the inverse image
of $(x_0,\ldots , x_p)$ under the map $\Ff _{p/}\rightarrow A_0 \times \ldots
\times A_0$. We do not make the assumption that $\Ff _{0/}$ is a discrete set,
and indeed for $x\in ob(A)$ the simplicial set $\Ff _{0/}(x)$ is exactly the
one which is considered to be parametrized by the object $x$. We have a map
of simplicial sets
$$
\Ff _{p/}(x_0,\ldots , x_p) \rightarrow A_{p/} (x_0,\ldots , x_p)\times \Ff
_{0/}(x_p)
$$
given by the projection $\Ff \rightarrow A$ and the structural map for $\Ff$
with respect to the arrow $0\rightarrow p$ in $\Delta$ corresponding to the
last vertex. The ``(contravariant) cartesian condition'' is that the above map
should be a weak equivalence of simplicial sets. Note that the ``covariant
cartesian condition'' would be the same but using the structural map to $\Ff
_{0/}(x_0)$ rather than to $\Ff _{0/}(x_p)$.

A cartesian family corresponds to a weak functor $A^o\rightarrow S$ in much the
same way as the Segal condition encodes the notion of weak category:
the action of the space of morphisms $A_{1/}(x,y)$ is given by the diagram
$$
\Ff _{0/}(y) \times A_{1/}(x,y)\stackrel{\cong}{\leftarrow}
\Ff _{1/}(x,y) \rightarrow \Ff _{0/}(x),
$$
the second morphism being the structural morphism for the map $0\rightarrow 1$
in $\Delta$ corresponding to the first vertex. The higher $\Ff _{p/}$ encode
homotopy-coherent associativity of this action.

In \cite{SimpsonAspects} the notion of {\em cartesian family} is defined by
saying that it is a precartesian family which satisfies a certain quasi-fibrant
condition. This quasi-fibrant condition (which is analogous to the
classical notion of quasi-fibration and is somewhat similar to Rezk's notion of
``sharp map'' \cite{Rezk2}) is designed to guarantee that cartesian
families over
Segal precats can be glued together. This glueing property ensures
representability of the associated functor of Segal precats, and allows us to
define a Segal category $S'$ with the property that a morphism
$A^o\rightarrow S'$ is exactly the same thing as a contravariant cartesian
family
over $A$. In \cite{SimpsonAspects} it is shown that there is a natural morphism
$S\rightarrow S'$, that this is an equivalence of Segal categories, and that
$S'$ is fibrant; thus $S'$ is a canonical fibrant replacement for $S$.  This
fact means that ``weak families'' of simplicial sets parametrized by a Segal
category $A$, i.e. weak functors $A^o\rightarrow S$, may be viewed as cartesian
families. The proofs in   \cite{SimpsonAspects} are given in the context of
$n$-categories but the same work in the Segal category context (or more
generally for Segal $n$-categories \cite{HirschowitzSimpson}).

In practice, there is no essential difference between the notion of
precartesian family and the notion of cartesian family. Generally speaking, the
natural constructions that one can make are precartesian but not cartesian;
then one should make a fibrant replacement (which is consequently quasi-fibrant)
to get a cartesian family. We will systematically ignore this point in the
remainder of the paper, and speak only of precartesian families but use the
terminology ``cartesian family''. The reader should note that in order to be
precise, one must make fibrant replacements sometimes. Since these are
essentially unique (i.e. unique up to coherent homotopy) this doesn't pose
any homotopy-coherence problems. Of course one should check that the previous
phrase is true $(\otimes )$.

\subnum{Segal categories of presheaves}

The fundamental construction underlying SGA 4 \cite{SGA4} is the
Yoneda embedding of a category into the category of presheaves over itself.
We have the same thing for Segal categories. For this section I should
acknowledge the suggestion of A. Hirschowitz who pointed out that it would be
interesting to look at the notion of representable functor in the context of
$n$-categories. 	And J. Tapia who pointed out to me that this was the
fundamental
thing in SGA 4; he is working on an altogether different generalization of it.

Let $S$ be the simplicial category of fibrant simplicial sets, and let $S'$ be
its replacement by an equivalent fibrant Segal category. If $A$ is any Segal
category, put
$$
\widehat{A}:= \underline{Hom}(A ^o, S').
$$
Recall that $A^o$ is the ``opposite'' Segal category, with the same objects as
$A$ and obtained by putting
$$
A^o_{p/}(x_0,\ldots , x_p):= A_{p/}(x_p,\ldots , x_0).
$$
The first step is that we would like to construct a natural morphism
$$
h_A:A\rightarrow \widehat{A}.
$$
In view of the definition of the internal $\underline{Hom}(A ^o, S')$ (see
above), constructing the morphism $h_A$ is equivalent to constructing the
``arrow
family''
$$
Arr _A: A^o \times A \rightarrow S'.
$$
We give two discussions of the construction of $Arr _A$. Both of these
constructions were done for $n$-categories in \cite{SimpsonAspects}.
We should also note that in the simplicial case, the ``arrow family'' is
certainly very classical; among other things it occurs in
Cordier-Porter \cite{CordierPorter}.

The easy case is when
$A$ is a strict simplicial category with fibrant simplicial $Hom$ sets. In this
case, the formula
$$
Arr _A(x,y) := A_{1/}(x,y)
$$
defines in an obvious way a morphism of strict simplicial categories
$$
A^o \times A \rightarrow S.
$$
There is a canonical fibrant replacement within the category of
simplicial sets, compatible with direct product (namely taking the singular
complex of the topological realization of a simplicial set), so we obtain a way
of replacing any simplicial category by one whose simplicial
$Hom$ sets are fibrant. This can be composed with the Dwyer-Kan strictification
described above, so if $A$ is any Segal category then we can replace $A$ by an
equivalent strict simplicial category with fibrant $Hom$ spaces and then apply
the construction of $Arr _A$ given in the present paragraph. Thus this
construction technically speaking suffices in order to define the morphism
$h_A$ and the reader wishing to avoid technicalities may skip the subsequent
paragraph.

The more complicated case is to treat directly the case where $A$ is a Segal
category. This has the advantage of avoiding a number of equivalences used in
the previous paragraph; however it makes use of the notion of ``cartesian
family'' described above (and for which the reader must refer to
\cite{SimpsonAspects}). We choose for fibrant replacement that $S'$
which was obtained using the notion of cartesian family. Thus, in order to
define the morphism
$$
Arr _A: A^o \times A \rightarrow S',
$$
we have to define a contravariant cartesian family over $A\times A^o$. We
do this
by first defining a natural precartesian family $\Ff$, then replacing by a
fibrant replacement $\Ff '$.   The precartesian family $\Ff$ has the very simple
formula $$
\Ff _{p/}((x_0, y_0),\ldots , (x_p, y_p)):=
A_{2p +1 /}(x_0,\ldots , x_p, y_p,\ldots , y_0).
$$
Note that
$$
(A\times A^o )_{p/}((x_0, y_0),\ldots , (x_p, y_p))=
A_{p/}(x_0,\ldots , x_p)\times A_{p/}(y_p,\ldots , y_0).
$$
The Segal condition for $A$ implies that the map
$$
\Ff _{p/}((x_0, y_0),\ldots , (x_p, y_p)) \rightarrow
A_{p/}(x_0,\ldots , x_p)\times A_{p/}(y_p,\ldots , y_0)
\times A_{1/}(x_p,y_p)
$$
is an equivalence. This is the cartesian condition for $\Ff$, so $\Ff$ is a
precartesian family. The morphism $Arr _A$ is defined by choosing a fibrant
replacement $\Ff '$ for $\Ff$.

In the above discussion, the Segal category $A$ must be small. For a
``big'' Segal category (by which we always mean one in which the objects can
form a class, but in which the $A_{p/}(x_0,\ldots , x_p)$ are still sets),
it doesn't seem to be reasonable to define $\widehat{A}$. However, we will run
across the following intermediate situation: suppose
$$
C\rightarrow A
$$
is a morphism from a small Segal category $C$ to a ``big'' Segal category $A$.
Then we still obtain a morphism
$$
i: A \rightarrow \widehat{C}.
$$
Define this by exhausting $A$ by small Segal categories $A_{\beta}$, and on each
of these define $i$ as the composition
$$
A_{\beta} \rightarrow \widehat{A}^{\beta} \rightarrow \widehat{C}.
$$

Here is the statement of our main ``Yoneda-type'' theorem.

\begin{theorem}
\label{yoneda}
If $A$ is any small Segal category then the morphism
$$
h_A: A\rightarrow \widehat{A}
$$
is fully faithful.
\end{theorem}
{\em Proof:}
We prove the following more general statement: if $G\in \widehat{A}$ and
if $x\in A$ then there is a natural equivalence
$$
\widehat{A}_{1/}(h_A(x), G) \cong G(x)
$$
(which is required to be compatible with $h_A$, see below).

We first point out how to go from here to the statement of the theorem:
for $x,y\in ob(A)$, apply the above to $G:= h_A(y)$.
We get
$$
\widehat{A}_{1/}(h_A(x), h_A(y))\cong h_A(y),
$$
but $h_A(y) \cong A_{1/}(x,y)$ by construction (recall that $h_A$ comes from
the arrow family). Thus
$$
\widehat{A}_{1/}(h_A(x), h_A(y))\cong A_{1/}(x,y).
$$
This equivalence will be compatible with the morphism $h_A: A\rightarrow
\widehat{A}$, so it shows that $h_A$ is fully faithful.

Now we show how to prove the more general statement. We can view $G$ as being a
cartesian family over $A$. In order to define a morphism
$$
G(x)\rightarrow \widehat{A}_{1/}(h_A(x), G)
$$
we need to define a morphism
$$
\Upsilon (G(x))\rightarrow \widehat{A}
$$
or equivalently a morphism
$$
[\Upsilon (G(x))\times A]^o \rightarrow S'
$$
restricting over $0\times A^o$ to $G$, and restricting over
$0\times A^o$ to $G$. This latter morphism corresponds to  a
contravariant cartesian family
$$
\Ff \rightarrow \Upsilon (G(x))\times A,
$$
with $\Ff$ restricting as above to $Arr _A(-,x)$ and $G$ on the endpoints.
In order to define the family $\Ff$, given that we already know its
restrictions to $0\times A$ and $1\times A$, it suffices to define
$$
\Ff _{p/} (u_0,\ldots , u_a; v_0,\ldots , v_b):= G_{p+1/}(u_0,\ldots , u_a,
v_0, \ldots , v_b, x).
$$
for $a,b\geq 0$ and $a+b+1=p$. Here $u_i, v_j \in ob(A)$ and the variables
$u_i$ indicate objects considered in $0\times A$; the variables $v_j$ indicate
objects considered in $1\times A$.
The simplicial restriction maps are obtained by those of $G$ whenever the
sequence of objects still contains an object of $0\times A$, otherwise
it is obtained by composing with the morphism $G\rightarrow A$. The structural
morphism to $\Upsilon (G(x))\times A$ will be seen in the upcoming verification.
We check the cartesian condition:
$$
G_{p+1/}(u_0,\ldots , u_a,
v_0, \ldots , v_b, x)\cong
A_{p+1/}(u_0,\ldots , u_a,
v_0, \ldots , v_b, x)\times G(x)
$$
$$
\cong
A_{p/}(u_0,\ldots , u_a,
v_0, \ldots , v_b)\times h_A(x)(v_b)\times G(x)
$$
$$
\cong
[\Upsilon (G(x))\times A]_{p/}(u_0,\ldots , u_a,
v_0, \ldots , v_b)\times h_A(x)(v_b).
$$
Thus $\Ff$ is a precartesian
family. As said previously, we are ignoring the difference between cartesian
and precartesian families. Thus we have defined our morphism
$$
G(x)\rightarrow \widehat{A}_{1/}(h_A(x), G).
$$

The next step is to define a morphism in the other direction:
$$
\widehat{A}_{1/}(h_A(x), G)\rightarrow G(x).
$$
For this, note that the restriction along $\{ x\} \rightarrow A$ gives a
morphism
$$
\widehat{A} \rightarrow S'.
$$
We obtain a morphism
$$
\widehat{A}_{1/}(h_A(x), G)\rightarrow S'_{1/}(h_A(x)(x), G(x)).
$$
On the other hand, the identity element gives a morphism
$\ast \rightarrow h_A(x)(x)=A_{1/}(x,x)$, and ``composing'' with this
gives
$$
S'_{1/}(h_A(x)(x), G(x)) \rightarrow S'_{1/}(\ast , G(x))
\cong G(x).
$$
As usual this ``composition'' requires inverting some equivalences which come
up in the notion of Segal category. We don't write out the details of that here
(although this neglect doesn't actually merit a $\otimes$). We get
our morphism  $$
\widehat{A}_{1/}(h_A(x), G)\rightarrow G(x).
$$

To complete the proof, we have to say that these two morphisms are inverses up
to homotopy. In one direction it is basically easy (modulo struggling with the
details of the weak compositions everywhere) that the composition
$$
G(x) \rightarrow \widehat{A}_{1/}(h_A(x), G)\rightarrow G(x)
$$
is homotopic to the identity of $G(x)$. For this direction, one way to proceed
would be to note that, for an appropriate Dwyer-Kan-Smith strictification and
then Dwyer-Hirschhorn-Kan cofibrant replacement, $A$ can be assumed to be a
strict simplicial category and $G$ a strict diagram $A\rightarrow S$. In this
setup we obtain (by just simplicially-enriching the easy discussion for
$1$-categories) a sequence
$$
G(x) \rightarrow Hom (A, S) _{1/}(h_A(x), G) \rightarrow G(x)
$$
whose composition is the identity of $G(x)$ on-the-nose.
In this formula, the simplicial category $Hom (A, S)$ is not necessarily the
``right'' one but it maps into $\widehat{A}$, and this is sufficient to
check that the above composition that we are interested in, is homotopic to the
identity. Note that the morphism in the strictified setup is homotopic to the
morphism we have constructed in the original weak situation.

It is somewhat more problematic to see why the composition
$$
\widehat{A}_{1/}(h_A(x), G)\rightarrow G(x)
\rightarrow \widehat{A}_{1/}(h_A(x), G)
$$
is the identity. This is because it is not clear (to me at least) whether
all of $\widehat{A}_{1/}(h_A(x), G)$ can in some way---and after appropriate
replacements of $A$ and $G$---be supposed to consist entirely of strict natural
transformations between strict diagrams.

Instead, we again make a more general statement, namely the naturality
of the morphism
$$
G(x) \rightarrow \widehat{A}_{1/}(h_A(x), G)
$$
in the variable $G$. This says that if $F$ and $G$ are objects
in $\widehat{A}$ then the diagram
$$
\begin{array}{ccc}
F(x) \times \widehat{A}_{1/}(F,G) & \rightarrow & G(x) \\
\downarrow && \downarrow \\
\widehat{A}_{1/}(h_A(x), F) \times \widehat{A}_{1/}(F,G)
& \rightarrow & \widehat{A}_{1/}(h_A(x), G)
\end{array}
$$
commutes up to homotopy.

For this statement and its proof, we first take note of the following remark:
if $U$ and $V$ are diagrams in $\widehat{A}$ then a morphism
$E\rightarrow \widehat{A}_{1/}(U,V)$ is by definition a morphism
$$
\Upsilon (E) \rightarrow \underline{Hom}(A^o, S')
$$
or equivalently a contravariant cartesian family over
$$
A\times \Upsilon (E)
$$
restricting to $V$ on $A\times 0$ and to $U$ on $A\times 1$
(in this last reduction we use the natural isomorphism $\Upsilon (E)^o\cong
\Upsilon (E)$ which interchanges $0$ and $1$). It is easy to see that a
precartesian family over $A\times \Upsilon (E)$, with restrictions $U$ and $V$,
is exactly the same thing as a precartesian family over $A\times I$ with
restrictions $V$ on $A\times 0$, and $U\times E$ on $A\times 1$.

With this remark in mind, we can return to the above diagram and (applying
the remark to the vertical arrows) note that it is the same thing as giving a
diagram in $\widehat{A}$ of the form
$$
\begin{array}{ccc}
h_A(x)\times F(x) \times \widehat{A}_{1/}(F,G) & \rightarrow & h_A(x)\times
G(x) \\
\downarrow && \downarrow \\
F \times \widehat{A}_{1/}(F,G)
& \rightarrow & G .
\end{array}
$$
Again applying the remark of the previous paragraph (but to the horizontal
arrows this time ) with $E= \widehat{A}_{1/}(F,G)$
we get that the above diagram is the same thing as a precartesian family
over
$$
A\times I \times \Upsilon (E)
$$
whose restrictions to the corners are respectively:
$$
A\times 0 \times 0: \;\;\; G
$$
$$
A\times 0 \times 1: \;\;\; F
$$
$$
A\times 1 \times 0: \;\;\; h_A(x)\times G(x)
$$
$$
A\times 1 \times 1: \;\;\; h_A(x)\times F(x).
$$
The restrictions to the edges $A\times I \times 0$ and $A\times I \times 1$
should be the cartesian families constructed above for $G$ and $F$ respectively;
the restrictions to $A\times 0 \times \Upsilon (E)$ and
$A\times 1 \times \Upsilon (E)$ should be the tautological families. The
construction of this cartesian family is done in the same way as the previous
construction for the morphism $h_A(x)\times G(x)\rightarrow G$, but starting
with the tautological cartesian family over $A\times \Upsilon (E)$
corresponding to the morphism $F\times E\rightarrow G$. We leave it to the
reader to write down the details $(\otimes )$. This gives the
homotopy-commutative diagram of naturality.

Let's now look at how to go from the above naturality statement to the fact
that our composition of morphisms is homotopic to the identity. For this,
apply the naturality statement with $F= h_A(x)$ and $G$ as given.
Then the naturality statement is a diagram
$$
\begin{array}{ccc}
A_{1/}(x,x) \times \widehat{A}_{1/}(h_A(x),G) & \rightarrow & G(x) \\
\downarrow && \downarrow \\
\widehat{A}_{1/}(h_A(x), h_A(x)) \times \widehat{A}_{1/}(h_A(x),G)
& \rightarrow & \widehat{A}_{1/}(h_A(x), G) .
\end{array}
$$
Plugging in the identity from $x$ to $x$  we get a
diagram
$$
\begin{array}{ccc}
\widehat{A}_{1/}(h_A(x),G) & \rightarrow & G(x) \\
\downarrow && \downarrow \\
\{ 1_{h_A(x)} \} \times \widehat{A}_{1/}(h_A(x),G)
& \rightarrow & \widehat{A}_{1/}(h_A(x), G) .
\end{array}
$$
The composition along the top followed by the right is the morphism we are
interested in; the other composition is the identity. Therefore,
homotopy-commutativity of the square shows that the composition
in question
$$
\widehat{A}_{1/}(h_A(x),G) \rightarrow  G(x)
\rightarrow \widehat{A}_{1/}(h_A(x),G)
$$
is the identity. This completes the proof of the theorem.
\eop

A lemma which will be used in several places below (and indeed, which is
at the origin of the statement of Theorem \ref{main}) is the following
calculation of $\widehat{C}$. Recall that the Heller model category of
simplicial presheaves \cite{Heller} was the precursor to the now standard
Joyal-Jardine model category \cite{Joyal} \cite{Jardine}; Heller's result
was the
special case of a category with trivial Grothendieck topology (which is the case
we need here).

\begin{lemma}
\label{heller}
Suppose $C$ is a small $1$-category. Then $\widehat{C}$ is equivalent to
$L(M)$ where $M = S^C$ is the Heller model category of simplicial presheaves
over $C$.
\end{lemma}
{\em Proof:} This is a special case of Th\'eor\`eme 12.1 of
\cite{HirschowitzSimpson}. To obtain the special case, replace $n$ by $0$ in
the statement of that theorem, and note that (in the notation of
\cite{HirschowitzSimpson}) a ``Segal $0$-category'' is the same thing as a
simplicial set.  One should also refer to Th\'eor\`eme 11.11 of the same
reference.

This statement is also given by Rezk in \cite{Rezk}, and the proof Rezk gives
uses some results of Dwyer-Kan. (The results of Dwyer-Kan were of course
much prior to \cite{HirschowitzSimpson}).
\eop

\subnum{Adjoint functors}

There is a notion of adjunction between functors of simplicial categories
or Segal categories, which is a direct generalization of the classical notion of
adjunction of functors. In making this generalization, it is best to specify
only one of the adjunction transformations and impose the condition that it
induces an equivalence between the appropriate simplicial $Hom$ sets. If one
tried to specify both of the classical adjunction transformations, this would
run into the homotopy-coherence problem that it would be necessary (in order
to obtain a well-behaved notion) to specify higher order homotopy coherencies.

The basic historical reference for this section is Cordier and Porter
\cite{CordierPorter}, who treat the case of adjunctions of homotopy-coherent
functors between simplicial categories. This should be completely equivalent to
what we say here. Furthermore, refering to their approach might allow easy
removal of the many $\otimes$ which appear in the following discussion.

Suppose $A,B$ are Segal categories (which we suppose fibrant) and
suppose $F: A\rightarrow B$ and $G: B\rightarrow A$ are functors. Suppose
$\eta : 1_B \rightarrow FG$ is a natural transformation; technically speaking,
this means
$$
\eta \in \underline{Hom}(B,B) _{1/}(1_B, FG),
$$
which in turn means that $\eta$ is a morphism of Segal categories
$$
B\times I \rightarrow B
$$
restricting to $1_B$ on $B\times 0$ and to $FG$ on $B\times 1$. Here as
throughout, $I$ denotes the category with two objects $0,1$ and a single
(non-identity) morphism $0\rightarrow 1$. Generally we consider $I$ as a Segal
category.

We obtain the following morphisms:
$$
\underline{Hom}(A^o\times A, S')
\stackrel{(G^o\times 1)^{\ast}}{\rightarrow}
\underline{Hom}(B^o \times A,S'),
$$
and
$$
\underline{Hom}(B^o\times B, S')
\stackrel{(1\times F)^{\ast}}{\rightarrow}
\underline{Hom}(B^o \times A,S').
$$
In particular, we have two elements
$$
(G^o\times 1)^{\ast}(Arr_A), \;\; (1\times F)^{\ast}(Arr _B)
\; \in
\underline{Hom}(B^o \times A,S').
$$
These represent respectively
$$
(x,y)\mapsto A_{1/}(Gx, y)
$$
and
$$
(x,y)\mapsto B_{1/}(x, Fy).
$$
In the same way as for the classical $1$-category case, the natural
transformation $\eta$ gives rise to a morphism ${\bf adj}(\eta )$ in the Segal
category $\underline{Hom}(B^o \times A,S')$ relating the above two
elements; this
morphism arises as a morphism $I\times B^o \times A\rightarrow S'$
restricting to $(G^o\times 1)^{\ast}(Arr_A)$ over $0\in I$ and
to $(1\times F)^{\ast}(Arr _B)$ over $1\in I$.

In the case where $F$ and $G$ are strict morphisms of strict simplicial
categories and $\eta$ is a strict natural transformation between them, the
adjunction morphism ${\bf adj}(\eta )$ is easy to describe; it is just
given by exactly the same formula as in the classical case.

The paragraph which follows contains a more technical description of how to
construct the morphism refered to above, in our framework of Segal categories.
This construction in turn relies on the explicit construction of a certain
cartesian family $\Phi _F$ which is left to the intrepid reader. The less
intrepid who are willing to accept that everything works as usual, may skip the
following paragraph.

Note that $(B\times I)^o \cong B^o \times I$ using $I^o \cong I$ (an involution
which switches $0$ and $1$).
Look in $\underline{Hom}(I \times B^o \times A, S')$
at
$$
(\eta \times F )^{\ast}(Arr _B).
$$
Over $0\in I$ this restricts to $(1_B\times F)^{\ast}(Arr _B)$.
Over $1$ this restricts to
$$
((FG)^o\times F)^{\ast}(Arr _B).
$$
Essentially speaking, this means that we have a natural transformation
$$
B_{1/}(FGx, Fy) \rightarrow B_{1/}(x, Fy).
$$
Note that $(FG)^o \times F$ is the composition
$$
B^o \times A
\stackrel{G^o \times 1}{\rightarrow}
A^o \times A \stackrel{F^o \times F }{\rightarrow}
B^o \times B .
$$
Thus
$$
((FG)^o\times F)^{\ast}(Arr _B)= (G^o\times 1)^{\ast} ((F^o \times F)^{\ast}Arr
_B).
$$
The morphism of functoriality for $F$ is a natural transformation
$$
A_{1/}(x,z) \rightarrow B_{1/}(Fx, Fz),
$$
which translates in our language to a morphism
$$
\Phi _F: I \times A^o \times A  \rightarrow S',
$$
restricting over $0$ to $Arr _A$, and
over $1$ to $(F^o \times F)^{\ast}Arr_B$.
Technically speaking, $\Phi _F$ needs to be constructed as a cartesian family
(recall that $Arr _A$ and $Arr _B$ are themselves cartesian families).
We leave this construction to the reader $(\otimes )$.
Now look at
$$
(1\times G^o\times 1)^{\ast}(\Phi _F): I \times B^o \times A\rightarrow S'.
$$
Heuristically it is the natural transformation
$$
A_{1/}(Gx,y) \rightarrow B_{1/}(FGx, Fy).
$$
We can ``compose'' this with the previous transformation to obtain a
natural transformation
$$
A_{1/}(Gx,y) \rightarrow B_{1/}(FGx, Fy)\rightarrow B_{1/}(x, Fy).
$$
Technically speaking, this means using the above two morphisms to
give the $01$ and $12$ edges which can be filled in to a morphism
$$
\Upsilon ^2(\ast , \ast ) \times B^o \times A\rightarrow S',
$$
the third ($02$) edge of which is a morphism
$$
{\bf adj}(\eta ): I\times B^o \times A\rightarrow S'
$$
restricting on the endpoints to $(1\times G)^{\ast}(Arr _A)$
and $(F^o \times 1)^{\ast}(Arr _B)$ respectively. This is the technical
description of how we get from the natural transformation
$\eta : 1_B\rightarrow FG$ to a natural
transformation
$$
{\bf adj}(\eta )(x,y): A_{1/}(Gx,y) \rightarrow B_{1/}(x, Fy).
$$

Now getting back to our discussion of adjoint functors, we say that {\em $\eta$
is an adjunction between $F$ and $G$} if the natural transformation
${\bf adj}(\eta )$ is an equivalence between
$(G^o \times 1)^{\ast}(Arr _A)$
and $(1 \times F)^{\ast}(Arr _B)$
(by ``equivalence'' here we mean internal equivalence in the Segal category
$\underline{Hom}(B^o \times A, S')$).

{\bf Remark:} In order to check the adjunction condition, it suffices to check
that for every pair of objects $x\in ob(B)$ and $y\in ob(A)$, the morphism
$$
{\bf adj}(\eta )(x,y): A_{1/}(Gx,y) \rightarrow B_{1/}(x, Fy)
$$
is a weak equivalence of simplicial sets. This is a general fact about
natural transformations between functors of Segal categories: being a levelwise
equivalence implies being an equivalence. It is Corollary 2.5.8 of
\cite{SimpsonLimits} (which was stated for $n$-categories but which
works the same way for Segal categories); a similar early result was shown in
\cite{SimpsonFlexible}.

\begin{lemma}
Suppose $F: A\rightarrow B$, $G: B\rightarrow A$ are functors of fibrant
Segal categories, and $\eta :B\times I \rightarrow B$ is a natural
transformation $1_B\rightarrow FG$ which is an
adjunction. Suppose that $C$ is another Segal category. Let $F_C, G_C$ be the
induced functors between $\underline{Hom}(C, A)$ and
$\underline{Hom}(C, B)$, and let $\eta _C$ denote the functor
$$
\underline{Hom}(C, B) \times I \rightarrow \underline{Hom}(C,B)
$$
defined by the composed morphism
$$
\underline{Hom}(C, B) \times I \times C
=
C\times \underline{Hom}(C, B) \times I
\rightarrow B\times I \stackrel{\eta}{\rightarrow} B.
$$
Then $\eta _C$ is a natural transformation
$$
1_{\underline{Hom}(C,B)}\rightarrow F_C G_C,
$$
which is an adjunction between $F_C$ and $G_C$.
\end{lemma}
{\em Proof:}
After the details of how to define everything, we will end up
with a natural transformation
$$
{\bf adj}(\eta _C)(u,v): \underline{Hom}(C,A)_{1/}(Gu,v) \rightarrow
\underline{Hom}(C,B)_{1/}(u, Fv).
$$
According to the previous remark, we have to show that this is an equivalence
for every $u: C\rightarrow B$ and $v: C\rightarrow A$. To check this, note that
$$
\underline{Hom}(C,A)_{1/}(Gu,v)
$$
is calculated by a homotopy-coherence calculation using the
$$
A_{1/}(Gu(c), v(c'))
$$
for $c,c'$ in $C$ (something like a coend, see Cordier-Porter
\cite{CordierPorter}). Similarly,
$$
\underline{Hom}(C,B)_{1/}(u,Fv)
$$
is calculated by the samehomotopy-coherence calculation using
$$
B_{1/}(u(c), Fv(c')).
$$
The fact that the adjunction induces an equivalence
$$
A_{1/}(Gu(c), v(c'))\cong B_{1/}(u(c), Fv(c'))
$$
for any $c,c'\in ob(C)$, implies that the two calculations give the same
answer; thus
${\bf adj}(\eta _C)(u,v)$ is an equivalence. This completes the proof, but a
number of details need to be followed through $(\otimes )$.
\eop

{\bf Construction:}
We can apply this to the case where $C=A$, where $u=F$ and where $v=1_A$.
We obtain an equivalence
$$
{\bf adj}(\eta _A)(F,1_A):\underline{Hom}(A,A)_{1/}(GF,1_A)
\stackrel{\cong}{\rightarrow} \underline{Hom}(A,B)_{1/}(F,F).
$$
In particular, there is an essentially unique element
$$
\zeta \in \underline{Hom}(A,A)_{1/}(GF,1_A)
$$
which goes to $1_F$ under the above equivalence. (To be more precise, what is
essentially unique---i.e. parametrized by a contractible space---is the pair
consisting of $\zeta$ plus a path between the image of $\zeta$ and $1_F$).

We leave it to the reader $(\otimes )$ to check that $\zeta$ is an adjunction
morphism going in the other direction between $F$ and $G$ (reversing the
appropriate things in the above discussion/definition).  We will use this
construction of the other adjunction morphism, at some point in the argument
below.

\begin{lemma}
\label{adjunctioncompositions}
With the above notations, the composed morphisms
$$
F \stackrel{\eta _{F()}}{\rightarrow} FGF \stackrel{F(\zeta )}{\rightarrow} F
$$
and
$$
G\stackrel{G(\eta )}{\rightarrow}GFG \stackrel{\zeta _{G()}}{\rightarrow} G
$$
are homotopic to the identity natural transformations of $F$ and $G$
respectively.
\end{lemma}

We don't give a proof of this here $(\otimes )$.

For the above places where details are left out in our discussion of
adjunction, the necessary arguments can probably be obtained from
Cordier-Porter \cite{CordierPorter}.

\subnum{Homotopy colimits}

It would be impossible to give a complete list of references to everything
pertaining to homotopy colimits (and limits). A non-exhaustive list includes
\cite{BousfieldKan} \cite{Vogt} \cite{Vogt2} \cite{EdwardsHastings}
\cite{Hirschhorn} \cite{DHK} \ldots .

Recall the notion of {\em homotopy colimit} in a simplicial category or Segal
category. If $A$ is a Segal category (which we suppose fibrant) and if $J$ is
a small Segal category, then we can form the {\em category of diagrams}
$\underline{Hom}(J,A)$. This is the ``homotopically correct'' one if $A$ is
fibrant. There is a morphism $c_J: A\rightarrow \underline{Hom}(J,A)$
induced by the projection $J\rightarrow \ast$; thus $c_J(x)$ is
the constant diagram with values $x$. Suppose $F: J\rightarrow A$ is a diagram.
If $x$ is an object of $A$ and $f: F \rightarrow c_J(x)$ is a morphism,
then we say that {\em $x$ is the homotopy colimit of the diagram $F$} and write
$$
(x,f)= colim _J(F)
$$
(or just $x= colim _J(F)$ if there is no confusion about $f$), if
for any object $y$ of $A$, the morphism of ``composition with $f$'',
which can be seen as the composition
$$
A_{1/}(x,y) \rightarrow \underline{Hom}(J, A)_{1/}(c_J(x), c_J(y))
\rightarrow
\underline{Hom}(J, A)_{1/}(F, c_J(y)),
$$
is an equivalence of simplicial sets. Here the second morphism is essentially
well-defined as ``composition'' in the Segal category
$\underline{Hom}(J, A)$, see above.

Note that we never speak of actual limits or colimits in a simplicial category,
so the notation $colim$ means homotopy colimit. If we forget to include the
qualifier ``homotopy'' in front of the word ``colimit'' in the text below,
the reader will insert it.  However, for homotopy limits or colimits of
simplicial sets, we keep the classical notation $holim$ or $hocolim$ so as not
to confuse these with $1$-limits or $1$-colimits in the $1$-category of
simplicial sets.

Note that
$$
\underline{Hom}(J, A) _{1/}(F, c_J(y))
\cong holim _{j\in J}A_{1/}(F(j), y)
$$
where the $holim$ on the right is the homotopy limit of simplicial sets.
Thus, we can rewrite the condition for being a homotopy colimit as saying
that for any object $y$, the composition morphism with $f$ gives an equivalence
of simplicial sets
$$
A_{1/}(x,y) \stackrel{\cong}{\rightarrow} holim _{j\in J}A_{1/}(F(j), y).
$$
In this sense, the homotopy colimit is in a certain sense dual to the homotopy
limit on the level of the simplicial $Hom$-sets of $A$ (i.e. the $A_{1/}(-,-)$).
In particular, we can verify certain formulae for homotopy colimits by
verifying the dual formulae for homotopy limits of simplicial sets. For example
it follows from \cite{Vogt2} that homotopy colimits commute with
other homotopy colimits.

There is an analogous definition of homotopy limit which we leave to the reader
to write down in our current language.

We now remark that colimits over a Segal category $J$ can be transformed into
colimits over a $1$-category $J'$; thus, in the above discussion, there would
be no loss of generality in considering the indexing category $J$ to be a
$1$-category. This remark follows from the following statement,
which we isolate as a lemma because it will also be used in the proof of the
main theorem.

\begin{lemma}
\label{transform}
If $C$ is a Segal category (which we assume fibrant), then there is a
$1$-category $D$ and a morphism $D\rightarrow C$ such that for any fibrant Segal
category $A$, the induced morphism
$$
\underline{Hom}(C,A) \rightarrow \underline{Hom}(D,A)
$$
is fully faithful. Furthermore, we can assume that $D$ is a ``Reedy poset'',
i.e. a poset with a Reedy structure such that the Reedy function is compatible
with the ordering.
\end{lemma}
{\em Proof:} It suffices to construct a
$1$-category $D$ and a subcategory $W\subset D$ with a morphism $D\rightarrow C$
sending the arrows of $W$ to equivalences in  $C$, such that this morphism
induces an equivalence
$$
L(D,W) \stackrel{\cong}{\rightarrow} C.
$$
To see that this suffices, recall from (\cite{HirschowitzSimpson}
Proposition 8.6--Corollaire 8.9) which in turn comes from (\cite{SimpsonLimits}
Theorem 2.5.1), that
$$
\underline{Hom}(L(D, W), A) \subset
\underline{Hom}(D, A)
$$
is the saturated full Segal subcategory consisting of the morphisms
$D\rightarrow 	A$ which send the morphisms of $W$ to equivalences in $A$.
(In the case $A= S'$ this result is  basically the same as the result of
Dwyer-Kan in \cite{DwyerKanDiags}).

Now for the construction of $D$ and $W$, we refer to \cite{HirschowitzSimpson}
Lemmes 16.1, 16.2. These basically say that one can construct $D$ and $W$ using
barycentric subdivision and the Grothendieck construction in the style of
Thomason.
\eop

{\em Caution:} One must be careful in combining this lemma with the Yoneda
result of Theorem \ref{yoneda}. In effect, one obtains (in the situation of the
lemma with $A=S'$) a sequence of three morphisms
$$
D\rightarrow C \rightarrow \widehat{C} \rightarrow \widehat{D}.
$$
The last two morphisms are fully faithful. The Yoneda morphism
$D\rightarrow \widehat{D}$ is also fully faithful. However, the composition of
these three morphisms is not in general the Yoneda morphism of $D$, so one
cannot
conclude that $D\rightarrow C$ must be fully faithful (which visibly it
isn't, in general). In fact, the composition of the above three morphisms is
homotopic to the Yoneda morphism for $D$ if and only if the original morphism
$D\rightarrow C$ is fully faithful.

\begin{corollary}
\label{transformcolim}
If $J$ is a Segal category (which we may assume fibrant) and if
$F: J\rightarrow A$ is a morphism to another Segal category, then
there is  a strict $1$-category (which we may assume to be a Reedy poset)
$J'$ and a morphism $g: J'\rightarrow F$ such  that if $colim ^A_{J'} F\circ g$
exists then $colim ^A_JF$ exists and the two colimits are equivalent.
\end{corollary}
{\em Proof:}
Choose $g: J'\rightarrow J$ (with $J'$ a Reedy poset) so that
$$
\underline{Hom}(J, S') \rightarrow \underline{Hom}(J', S')
$$
is fully faithful. Now note that if $G: J\rightarrow S'$ is a
simplicial set diagram over $J$, we have
$$
holim ^{S'} _J G \cong \underline{Hom}(J, S')_{1/}(\ast , G).
$$
The same holds for $J'$.
Therefore the fully faithful property implies that
$$
holim ^{S'} _J G \stackrel{\cong}{\rightarrow}
holim ^{S'} _{J'} G\circ g.
$$
Now the fact that homotopy colimits in $A$ are dual to homotopy limits of
the simplicial $Hom$ sets, implies that for any diagram $F:J\rightarrow A$,
$$
colim ^A_{J'}F\circ g \rightarrow
colim ^A_{J}F
$$
is an equivalence, and in fact existence of the first colimit implies existence
of the second one. For the statement about existence we use the fully
faithful property of the lemma (for target $A$ this time) to say that
$$
\underline{Hom}(J, A) _{1/}(F, c_J(colim ^A_{J'}F\circ g))
\rightarrow
\underline{Hom}(J', A) _{1/}(F\circ g, c_{J'}(colim ^A_{J'}F\circ g ))
$$
is an equivalence, so there exists a morphism of $J$-diagrams from
$F$ to $colim ^A_{J'}F\circ g$ restricting to the colimit morphism over $J'$;
now we can applly the previous discussion about $holim ^{S'}$ to get that this
morphism is a $J$-colimit.
\eop

We say that $A$ {\em admits all small homotopy colimits} if for any small Segal
category $J$ and for any diagram
$J\rightarrow F$, the homotopy colimit exists. From the previous corollary, it
suffices to check the existence of colimits over $1$-categories $J$ which we
can furthermore assume are Reedy posets.

\begin{lemma}
\label{samecolims}
Suppose $A\rightarrow B$ is a fully faithful morphism of Segal
categories. Suppose $F:J\rightarrow A$ is a diagram. If
$$
colim ^{B}_J(F)
$$
is in $A$, then the natural morphism
$$
colim ^{B}_J(F) \rightarrow
colim ^{A}_J(F)
$$
is an equivalence.
\end{lemma}
{\em Proof:}
The facts that $colim ^{B}_J(F)$ is in $A$ and that the inclusion of $A$
in $B$ is fully faithful imply that we have a morphism of $J$-diagrams
in $A$
$$
F \rightarrow c_J[colim ^{B}_J(F) ].
$$
Therefore there is up to homotopy a unique morphism
$$
colim ^A_J(F)\rightarrow colim ^{B}_J(F)
$$
whose composition with the canonical morphism of diagrams for the
$colim ^A_J$, is the above
morphism.   In the other direction, we have a morphism of
$B$-diagrams
$$
F\rightarrow c_J[colim ^A_J(F)].
$$
Again we get an essentially unique morphism
$$
colim ^{B}_J(F) \rightarrow colim
^A_J(F)
$$
(which is the morphism in the statement of the lemma). Essential uniqueness
implies that the compositions in both directions are homotopic to the identity,
thus our morphism is an equivalence.
\eop

\begin{lemma}
\label{objectbyobject}
If $C$ is a small Segal category, then homotopy colimits in $\widehat{C}$
exist and are calculated object-by-object.
\end{lemma}
{\em Proof:}
According to Lemma \ref{transform}, there is a small $1$-category
$D$ and a morphism $D\rightarrow C$ such that this induces a fully faithful
morphism
$$
\widehat{C}\rightarrow \widehat{D}.
$$
Furthermore, from the proof of Lemma \ref{transform}, we may assume that there
is a subcategory $W\subset D$ such that $C$ is equivalent to the
localization $L(D,W)$. This implies that $\widehat{C}$ is the full subcategory
of $\widehat{D}$ consisting of diagrams $X: D^o\rightarrow S'$ such that for
any arrow $w\in W$, $X(w)$ is an equivalence. This situation is identical to
that of Dwyer-Kan in \cite{DwyerKanDiags}, and all of the elements going into
here are due to \cite{DwyerKanDiags} in this case.

Next we recall that $\widehat{D}$ is equivalent to $L(M)$ where $M$ is the
Heller closed model category $S^D$ of simplicial presheaves over $D$
(Lemma \ref{heller}). Now, suppose we have a
diagram $F: J\rightarrow \widehat{C}$, which we may also consider as a diagram
in $\widehat{D}$. We may assume that $J$ is a Reedy poset. The argument of
\cite{HirschowitzSimpson} (see chapter 18 for example) allows us to
``strictify''
and assume that $F$ is the projection of a diagram $F': J\rightarrow M$.
Furthermore we may assume that $F'$ is Reedy-cofibrant in the variable $J$.
Then
$$
colim ^{\widehat{D}}_J(F) =colim ^{L(M)} _J(F)
$$
exists and is calculated by taking the $1$-colimit of $F'$ in $M$ (see the
discussion at the proof of $(i)\Rightarrow (ii)$ in Theorem \ref{main} below).
This $1$-colimit is calculated object-by-object over $D$ (recall that $M$ is the
category of simplicial presheaves on $D$). On the other hand, the Reedy
cofibrant
condition for $F'$ also holds object-by-object. Therefore for any $d\in ob(D)$,
the $1$-colimit of $F' (j)(d)$ over $j\in J$, is also the homotopy colimit. This
shows that the homotopy colimit is calculated object-by-object, i.e.
$$
colim ^{\widehat{D}}_J(F)(d) = colim ^{S'}_{j\in J}F(j)(d)
= hocolim _{j\in J}F(j)(d)
$$
for $d\in ob(D)$.
On the other hand, the
fact that $F$ is a diagram in $\widehat{C}$ means that for arrows $w$ in $W$,
and for any $j\in J$, we have that  $F'(j)(w)$ is an equivalence of simplicial
sets. Homotopy-invariance of the $1$-colimit of a Reedy cofibrant diagram
\cite{Hirschhorn} implies that the arrow
$$
colim ^{\widehat{D}}_J(F)(w)
$$
is an equivalence for any arrow $w$ in $W$. It follows that
$$
colim ^{\widehat{D}}_J(F)\in ob(\widehat{C}).
$$
Now by Lemma \ref{samecolims},
$$
colim ^{\widehat{C}}_J(F)\cong colim ^{\widehat{D}}_J(F)
$$
including the statement that the homotopy colimit in $\widehat{C}$ exists.
Finally, we have
$$
colim ^{\widehat{C}}_J(F)(d)
= colim ^{\widehat{D}}_J(F)(d) = colim ^{S'}_{j\in J}F(j)(d).
$$
This completes the proof.
\eop

\subnum{Smallness and rearrangement of colimits}

Recall from \cite{Hirschhorn} the notion of {\em sequential colimit}.
This is a colimit indexed by an ordinal $\beta$ (where the ordered set
$\beta$ is considered as a category with morphisms going in the increasing
direction) with the additional property that if $i\in \beta$ is a limit element
then the $i$-th object $X_i$ is equivalent to the colimit of the $X_j$ for
$j<i$.
A diagram giving rise to a sequential colimit will be called a {\em sequential
diagram}. In giving these definitions for a Segal category $A$, the notion of
colimit which occurs is the notion of homotopy colimit as defined above.

A diagram or colimit is {\em essentially sequential} if it satisfies the
sequential condition at sufficiently large points. In what follows we shall
make no distinction between sequential and essentially sequential (an
essentially sequential diagram can be replaced by a sequential one which gives
the same colimit, by starting out with a constant diagram in low degrees).

Suppose $A$ is a Segal category admitting all small colimits.An
object $z\in ob(A)$ is said to be {\em $\beta$-small in $A$} if for any ordinal
$\delta$ of size $\geq \beta$ and any sequential diagram $X: \delta \rightarrow
A$, the natural morphism
$$
hocolim _{i\in \delta} A_{1/}(z, X_i) \rightarrow A_{1/}(z, colim
^A_{\delta}X)
$$
is an equivalence. We say that $z$ is {\em small in $A$} if there is a cardinal
$\beta$ such that $z$ is $\beta$-small in $A$.

{\em Example:} Lemma \ref{objectbyobject} shows that the objects of $C$ are
small in $\widehat{C}$. On the other hand, every object of $\widehat{C}$ can
be expressed as a small homotopy colimit of objects of $C$ (see Lemma
\ref{express} below). From this  it follows easily that every object in
$\widehat{C}$ is small in $\widehat{C}$ (although of course there is no bound
uniform over the {\em class} $ob(C)$).

We will now treat some aspects of colimits which are useful in
connection with the notion of smallness.

Suppose $J$ is a small $1$-category, and $F: J\rightarrow A$ is a diagram.
Choose a well-ordering of the objects of $J$, in other words choose an ordinal
$\beta$ and an isomorphism $ob (\beta )\cong ob(J)$. We assume that $\beta$ is
the first ordinal of its cardinality. For each $i\in \beta$ let
$J_i$ denote the full subcategory of objects $<i$. Put
$$
X_i := colim _{J_i} (F|_{J_i}).
$$
Then the  $X_i$ form a sequential diagram, and we have
$$
colim _J (F) = colim _{\beta} (X_i).
$$
This can be proved by noting the dual property for homotopy limits of
simplicial sets.

We call the above expression a {\em normalized reindexing of the colimit}. The
word ``normalized'' refers to the condition that $\beta$ be the first ordinal of
its cardinality. With this condition, we get that each $X_i$ is a colimit of
size $< |\beta |$ (this latter notation is the cardinality of $\beta$).

Now we discuss another aspect of rearranging colimits. Let $A$ be a
Segal category admitting small colimits and let $C\subset A$ be a small full
subcategory. For any ordinal $\beta$ we will define a full subcategory
$$
A(< \beta C)\subset A,
$$
by the following prescription: it is the smallest saturated full subcategory of
$A$ containing $C$ and closed under colimits of size $< |\beta |$.

It is clear that if there is an ordinal $\beta ' < \beta$ of the same
cardinality, then $A(< \beta 'C) = A(< \beta C)$.

The following lemma is our main statement giving a normal form for
successive colimits.

\begin{lemma}
\label{rearrangement}
Suppose $\beta$ is the first ordinal of its cardinality. There are two cases.
\newline
(1) \, If $\beta$ is a limit of ordinals of strictly increasing cardinality,
then
$$
A(< \beta C) = \bigcup _{\gamma < \beta }  A(< \gamma C).
$$
(2) \, On the other hand, if $\beta$ is the limit of ordinals $\gamma$ all
having the same cardinality, then (letting $\gamma$ denote the first of these)
$A(< \beta C)$ is the saturated full subcategory of $A$ consisting of objects
which are $\gamma$-sequential colimits of objects of
$A(< \gamma C)$.
\end{lemma}
{\em Proof:}
In the first case (1), take the union in question, and note that it is indeed
closed under colimits of size $< \beta $, because any such colimit over $J$
has the property that there exists $\gamma < \beta $ with $|J| < \gamma$
so the colimit exists in $A(<\gamma C)$.

The main problem is to treat the second case. Let $A'\subset A$ be the
saturated full subcategory of $A$ consisting of objects which are
$\gamma$-sequential colimits of objects of $A(< \gamma C)$.
Suppose $F: J\rightarrow A'$ is a diagram of size $< \beta$. Note that
this implies that $|J|\leq |\gamma |$. Again we
isolate two cases: \newline
(a) where $|J| < |\gamma |$;
and \newline
(b) where $F$ is a sequential diagram taken over $J=\gamma$.

In these two cases, we will show that $colim _JF \in A'$. This suffices,
in view of the reorganization of the colimit over a $J$ of cardinality
$|\gamma |$.

In the first case, note that by Corollary \ref{transformcolim}, we can assume
that $J$ is a $1$-category which is a Reedy poset (note that the operation of
Corollary  \ref{transformcolim} doesn't
increase the size of $J$ beyond the countable cardinal).

Now by doing an induction on the Reedy function of the poset $J$ (and using the
assumption $|J| < |\gamma |$) we can rearrange the individual limits expressing
our objects $F(j)$ as objects in $A'$, so that  we have a doubly indexed system
$F_i(J)$ for $i \in \gamma $ and $j\in J$, such that
$$
F(j)= colim _{i\in \gamma} F_i(J).
$$
Here the systems in the variable $i$ are sequential. Now we may set
$G_i:= colim _J F_i(j)$. This colimit lies in $A(<\gamma C)$, and we have
$$
colim _JF = colim _{\gamma} G.
$$The
diagram $G$ may then be replaced by a sequential diagram giving the same
colimit. This proves that
$$
colim _JF \in A',
$$
so we have finished treating case (a).

For case (b), we can again (by induction on the ordered set $\gamma$)
suppose that our colimit comes from a doubly-indexed family
$F_i(j) \in A(<\gamma C)$ this time indexed by $\gamma \times \gamma$. By a
diagonal reindexing of this family we can express $colim _JF$ as
a $\gamma$-sequential colimit of objects of $A(<\gamma C)$, so again
$$
colim _JF \in A',
$$
and we have finished case (b). As remarked above, this suffices to obtain
case (2) of the lemma.
\eop

\subnum{Generating subcategories}

We have the following notions of generation: these are for a
saturated full sub-Segal category $C\subset A$, and we suppose that
$A$ is closed under colimits.

We say that {\em $C$ strongly generates $A$} if the morphism $A \rightarrow
\widehat{C}$ is fully faithful.

We say that {\em $C$ generates $A$ by colimits} if the smallest saturated full
sub-Segal  category of $A$ which contains $C$ and is closed under $A$-colimits,
is $A$ itself.

As a preliminary for the subsequent proposition we have the following lemma.
The main problem in giving the proof is to be careful to avoid
an error related to the ``caution'' after Lemma \ref{transform}. It should be
possible to give a more conceptual proof of this lemma using the ``coend''
construction of Cordier-Porter \cite{CordierPorter}.

\begin{lemma}
\label{express}
If $D$ is a small Segal category, then every object of $\widehat{D}$ can be
expressed as a small homotopy colimit of objects of $D$ (in the Yoneda
embedding).
\end{lemma}
{\em Proof:}
First suppose $D$ is a $1$-category. Then $\widehat{D}\cong L(M)$ where $M$ is
the Heller model category $S^D$ of simplicial presheaves over $D$ (Lemma
\ref{heller}). Any object of $M$ is equivalent to a simplicial object in the
category of formal disjoint unions of objects of $D$ (C. Teleman pointed this
out to me). From there we can go to an expression for any object as a homotopy
colimit of objects of $D$. This treats the case of a $1$-category.

Suppose now that $D$ is a Segal category, and choose a $1$-category $C$ with a
morphism $f: C\rightarrow D$ which induces a fully faithful morphism
$$
i: \widehat{D} \rightarrow \widehat{C}.
$$
We know the present lemma for $\widehat{C}$.
Plugging in $A:= \widehat{D}$ in the argument which will be given below in the
proof of Theorem \ref{main} (the part $(ii) \Rightarrow (iii)$) we obtain
existence of an adjoint $\psi : \widehat{C} \rightarrow \widehat{D}$.
Note that the present lemma is used in that argument, for $C$; but this we know
from the previous paragraph. Full-faithfulness of $i$ implies that
$\psi \circ i\cong 1_{\widehat{D}}$. Suppose now that $U\in \widehat{D}$.
Express
$$
iU \cong  colim ^{\widehat{C}}_Jh_C\circ F
$$
for a small diagram $F: J\rightarrow C$ composed with the Yoneda $h_C:
C\rightarrow \widehat{C}$. Now
$$
U\cong \psi i U \cong \psi (colim ^{\widehat{C}}_Jh_C\circ F).
$$
The fact that $\psi$ is an adjoint to $i$
implies that $\psi$ preserves colimits. Therefore we get
$$
U\cong colim ^{\widehat{D}}_J\psi \circ h_C\circ F.
$$
However, we have the formula
$$
\psi \circ h_C \cong h_D \circ f.
$$
(this follows from the adjunction property of $\psi$). Therefore we get
$$
U\cong colim ^{\widehat{D}}_Jh_D \circ f \circ F.
$$
In particular $U$ is a colimit of $h_D$ composed with the diagram
$f\circ F : D\rightarrow \widehat{D}$.
\eop

The following proposition gives an equivalence between our two notions of
generation.

\begin{proposition}
\label{generation}
Suppose $A$ is a Segal category closed under small homotopy colimits.
The following conditions are equivalent.
\newline
(a) there exists a small saturated full subcategory $C\subset A$ consisting of
objects which are small in $A$ and  which strongly generates $A$;
\newline
(b) there exists a small saturated full subcategory $C$ consisting of
objects which are small in $A$, such that $C$
generates $A$ by colimits.
\end{proposition}
{\em Proof:}
For $(a)\Rightarrow (b)$, suppose  $X\in ob(A)$ and  write
(by Lemma \ref{express})
$$
X = colim _J^{\widehat{C}}F
$$
for a diagram $F: J\rightarrow C$. The strong generation hypothesis that
$A\rightarrow  \widehat{C}$ is fully faithful implies (cf Lemma \ref{samecolims}
above) that $X$ is also the colimit  $colim ^A_J(F)$. Thus $C$ generates $A$ by
colimits (and in fact we obtain that every element of $A$ is expressed as a
single colimit of a diagram in $C$).
This shows that $(a) \Rightarrow (b)$.

We now show that $(b)\Rightarrow (a)$. Let $C$ be a saturated full
subcategory which generates $A$ by colimits, and suppose that there is a
cardinal $\gamma$ such that the objects of $C$ are $\gamma$-small in $A$. We
may also assume that the objects of $C$ are $\gamma$-small in $\widehat{C}$.
Recall that we have defined a saturated full subcategory $A(<\beta C)$, which is
the smallest one containing $C$ and closed under colimits of size $<\beta$.

We claim that for $\beta$ big enough (say bigger than $\gamma$), the
elements of $A(<\beta C)$ are $\beta$-small in $A$. Indeed, the objects of $C$
will be $\beta$-small, and it is easy to see that a colimit of size $<\beta$ of
objects which are $\beta$-small, is again $\beta$-small.

Now choose $\beta$ big enough, and set
$$
D:= A(<\beta C).
$$
From the previous paragraph, the elements of $D$ are $\beta$-small in $A$.
We claim that the morphism
$$
A \stackrel{i_D}{\rightarrow} \widehat{D}
$$
is fully faithful.
Suppose $Y$ is an object of $A$.
Look at the
saturated full subcategory $A'\subset A$ of objects $X$ such that
$$
A_{1/}(X,Y) \rightarrow \widehat{D} _{1/}(i_DX, i_DY)
$$
is an equivalence. We note first of all that $A'$ contains $D$.
To see this, note that if $X\in ob(D)$ then by definition
$$
i_D(Y) (X) = A_{1/}(X,Y).
$$
On the other hand, from Theorem \ref{yoneda}
$$
\widehat{D} _{1/}(i_DX, i_DY) \cong i_D(Y)(X).
$$
Therefore we obtain that $X\in A'$.

Next we claim that
$A'$ is closed under sequential colimits of size $\geq \beta$.
Suppose $X$ is a $\delta$-sequential $A$-colimit of $X_i$ with the
$X_i$ in $A'$, and with $|\delta |\geq \beta $.
We have
$$
A_{1/}(X, Y) = holim _{\delta} A_{1/}(X_i, Y).
$$
This maps  by an equivalence to
$$
lim _{\delta} \widehat{D}_{1/}(X_i , Y).
$$
In turn, this is equivalent to
$$
\widehat{D}_{1/}(colim _{\delta} ^{\widehat{D}}X_i , Y).
$$

Now we have a morphism in $\widehat{D}$
$$
colim _{\delta} ^{\widehat{D}}X_i\rightarrow X.
$$
We claim that this morphism is an equivalence. Indeed, suppose $Z$ is in $D$.
Then
$$
\widehat{D}_{1/}(Z, colim _{\delta} ^{\widehat{D}}X_i)
\cong
colim _{i\in \delta} ^S \widehat{D}_{1/}(Z,X_i)
$$
(filtered colimits of simplicial-set diagrams over $D$ are calculated
object-by-object). On the other hand,
$$
\widehat{D}_{1/}(Z,X) \cong X(Z) \cong A_{1/}(Z, X)
$$
and finally
$$
A_{1/}(Z,X) = A_{1/}(Z, colim ^A_{i\in \delta} X_i ) \cong
hocolim _{i\in \delta} A_{1/}(Z, X_i)
$$
the latter because $|\delta | \geq \beta $ and as remarked above,
the elements of $D$ are $\beta$-small in $A$. Finally note that
$$
A_{1/}(Z, X_i)=X_i(Z) \cong \widehat{D}_{1/}(Z,X_i).
$$
We obtain, after all of this, that the morphism
$$
colim _{\delta} ^{\widehat{D}}X_i\rightarrow X.
$$
induces an equivalence on morphism spaces from any object of $D$. This implies
that it is an equivalence, as claimed.

Recall from above that we have equivalences
$$
A_{1/}(X, Y) \cong lim _{\delta} A_{1/}(X_i, Y)
\cong
lim _{\delta} \widehat{D}_{1/}(X_i , Y).
\cong
\widehat{D}_{1/}(colim _{\delta} ^{\widehat{D}}X_i , Y).
$$
On the other hand the composed morphism factors as
$$
A_{1/}(X, Y)\rightarrow
\widehat{D}_{1/}(X,Y) \rightarrow
\widehat{D}_{1/}(colim _{\delta} ^{\widehat{D}}X_i , Y).
$$
From the claim of the previous paragraph, the second morphism is an
equivalence; therefore the first morphism is an equivalence, which shows that
$X\in A'$.

We have now shown that $A'$ is closed under sequential colimits.
The result of Lemma \ref{rearrangement} implies that $A'=A$. Therefore the
morphism $A\rightarrow \widehat{D}$ is fully faithful, and $D$ strongly
generates $A$ giving condition $(a)$ of the proposition. \eop

We now turn to generation for classes of morphisms.
Suppose $\Ff$ is a set of homotopy classes of morphisms in $A$. An {\em
$\Ff$-fibration} is a morphism $f:X\rightarrow Y$ which satisfies the weak
lifting property for elements of $\Ff$ (i.e. whose image satisfies the lifting
property in $ho(A)$). Here the morphism in $\Ff$ goes on the left in the square
diagram, and the morphism $f$ goes on the right. Suppose $\Ff$ is a subset of
homotopy classes of morphisms. The {\em class of morphisms generated (in terms
of lifting) by $\Ff$} is the largest subclass of (homotopy classes of) morphisms
$\overline{\Ff}$ such that the $\Ff$-fibrations are the same as the
$\overline{\Ff}$-fibrations.

\begin{lemma}
Suppose $\Ff$ is a set of homotopy classes of morphisms. Then any morphism
in $A$ which is a sequential limit of pushouts by morphisms in $\Ff$, is in
$\overline{\Ff}$.
\end{lemma}
{\em Proof:}
Lifting for morphisms in $\Ff$ implies lifting for morphisms which are pushouts
by morphisms in $\Ff$, and also for sequential colimits of such.
\eop

We will apply the above in the following case: where we are given a morphism
of Segal categories
$\psi : B \rightarrow A$. We will say that an arrow $f$ in $B$ (i.e. a vertex of
some $B_{1/}(x,y)$) is {\em $\psi$-trivial} if $\psi (f)$ is an internal
equivalence in $A$. We say that an arrow in $B$ is a {\em $\psi$-fibration} if
it is a fibration in the sense of the paragraph before the previous lemma, for
the class $\Ff$ of $\psi$-trivial morphisms. In the main theorem below, we
will be interested in when the class $\Ff$ of $\psi$-trivial morphisms is
generated (in terms of lifting as per the above definition) by a small subset of
$\psi$-trivial morphisms.

\subnum{The main theorem}

The following theorem is a generalization of Giraud's theorem
characterizing Grothendieck topoi.

\begin{theorem}
\label{main}
Suppose $A$ is a Segal category (which we may assume fibrant). The
following conditions are equivalent.
\newline
(i) \, There exists a cofibrantly generated closed model category $M$ such
that $A$ is equivalent to $L(M)$;
\newline
(ii)\, All  small homotopy colimits exist in $A$,
and there exists a cardinal $\beta$ and a small subset of objects $\Gg \subset
A_0$ such that the objects of $\Gg$ are $\beta$-small in $A$, and such that
$\Gg$
generates $A$ by colimits;
\newline
(iii)\, There exists a small $1$-category $C$ and a morphism $g:C\rightarrow A$
sending objects of $C$ to objects which are small in $A$,
which induces a fully faithful inclusion
$$
i:A\rightarrow \widehat{C};
$$
and there is a morphism $\psi : \widehat{C}\rightarrow A$ together with
a natural transformation
$$
\eta _X: X\rightarrow i\psi (X)
$$
such that $\eta$ induces an adjunction between $i$ and $\psi$.
\newline
(iv)  The category
$A$ admits all small homotopy colimits, and there exists a small
$1$-category $C$
and a functor $\psi : \widehat{C}\rightarrow A$ commuting with colimits,  such
that $A$ is the localization of $\widehat{C}$ by inverting the morphisms which
$\psi$ takes to equivalences, and such that the $\psi$-trivial morphisms of
$\widehat{C}$ are generated (in terms of lifting) by a small subset of
$\psi$-trivial morphisms of $\widehat{C}$.

Note that (i) implies that $A$ admits all small homotopy
limits too. In (iii), the fully faithful condition implies that the adjunction
morphism going in the other direction is an equivalence between $\psi \circ i$
and $1_A$.

We call a Segal category $A$ which satisfies these equivalent conditions,
an {\em
$\infty$-pretopos}.
If furthermore there exists $C\rightarrow A$ as in condition (iii)
such that the adjoint $\psi$ preserves finite homotopy limits, then we say that
$A$ is an {\em $\infty$-topos}.
\end{theorem}

The terminology ``$\infty$-topos'' naturally gives rise to a number of
conjectures, definitions, generalizations, predictions etc., which would be
too numerous to list here. As an example, we mention that
the applications of the theory of topoi to mathematical logic (cf e.g.
Moerdijk-MacLane \cite{MoerdijkMacLane}), should give rise to generalizations
in the case of $\infty$-topoi---which if they exist could be called
``$\infty$-categorical logic'' or ``higher-dimensional logic''.

{\em Remark:} It would be interesting to know what conditions on the closed
model category $M$ correspond to the $\infty$-topos condition.
(In a similar vein, it would be good to know that the $\infty$-topos condition
is independent of the choice of $C$ in condition (iii).) Charles Rezk
(\cite{RezkLetter} and later \cite{Rezk2}) points out that in the case of
presheaves on categories with Grothendieck pretopologies, the exactness of the
associated sheaf functor corresponds exactly to the condition that the
pretopology be a topology (this letter from Rezk to Hirschhorn was one of the
main elements motivating the present paper). We might expect a similar sort of
behavior here; one might even go so far as to conjecture that the $\infty$-topoi
are exactly the Segal categories of simplicial presheaves on Grothendieck sites.
Another conjecture (more reasonable) would be that the $\infty$-topoi correspond
exactly to the right proper closed model categories.

In keeping with the above remark, one should note that our theorem is not
strictly speaking an exact generalization of Giraud's theorem, because
we treat the case where $\psi$ may not be exact, and we obtain a weaker
result (existence of a closed model structure, rather than existence of  a
site).

P. Hirschhorn points out that if $M$ is a cofibrantly generated closed model
category, then the opposite category $M^o$ (which is again a closed model
category admitting all small limits and colimits) will not in general be
cofibrantly generated. Similarly, the opposite of an $\infty$-pretopos will
generally not be an $\infty$-pretopos. The point is that the generation
condition is asymmetric. For example (as Hirschhorn pointed out in an email) the
closed model category $Sets^o$ is not cofibrantly generated; indeed the only
small objects are those which correspond to the sets $\emptyset$ and $\ast$ (as
these are the only sets which are cosmall in $Sets$). These sets don't generate
$Sets$ by inverse limits, so they don't generate $Sets^o$ by colimits.

\bigskip

Our strategy for the proof of Theorem \ref{main} is similar to the strategy of
the proof of Giraud's theorem \cite{SGA4}: we prove
$$
(i)\Rightarrow (ii) \Rightarrow (iii) \Rightarrow (iv) \Rightarrow (i).
$$
For various reasons, our statements
(i)--(iv) are not precise generalizations of the statements which are numbered
in the same way in Giraud's theorem. Among other things, there is too much
``element-wise'' reasoning with respect to the $Hom$ sets in SGA 4 \cite{SGA4}
to make a direct generalization possible.
One major difference in our proof is that at the last step, we replace the
notion of Grothendieck topology by the notion of localization of a closed
model category. This explains why the theorem concerns
pretopoi rather than topoi.

\subnum{The proof of $(i) \Rightarrow (ii)$}

It is not difficult to see that if $M$ is a closed model
category (admitting small limits and colimits), then $L(M)$ admits all small
homotopy limits and homotopy colimits. A version of this statement was
known to Edwards and Hastings
\cite{EdwardsHastings}. Our technique comes out of the
Dwyer-Hirschhorn-Kan methods for calculating homotopy (co)limits and their
methods for calculating the function spaces in $L(M)$.
The argument we are about to present was done for
products and coproducts in Lemme 8.4 of \cite{HirschowitzSimpson}, but it
works the same way for arbitrary colimits (resp. limits).
We briefly review this
for the reader's convenience.

Suppose we want to
calculate the homotopy colimit of a diagram
$$
F: J\rightarrow L(M)'.
$$
It is shown in \cite{HirschowitzSimpson} that we may assume that
this diagram comes from a strict diagram $F: J\rightarrow M$. As in Corollary
\ref{transformcolim}, we may also assume that $J$ is a Reedy category and a
poset
with the ordering compatible with the Reedy structure. We may replace $F$ by a
levelwise equivalent Reedy-cofibrant diagram, so we may assume that $F$ is
Reedy-cofibrant. With these hypotheses we will show that the strict $1$-colimit
of $F$ in $M$,
$$
colim ^M_J(F)
$$
is a representative for the homotopy colimit $colim ^{L(M)}_J(F)$
which implies that the latter exists. To do this, consider an object $Y\in M$.
Choose a Reedy-fibrant simplicial resolution (see \cite{DHK}
or \cite{Hirschhorn})
$$
Y\rightarrow Z_{\cdot}.
$$
Recall from \cite{DHK} and \cite{Hirschhorn} that if $X$ is any object of $M$,
we have
$$
L(M)_{1/}(X,Y) \cong M_{1/}(X, Z_{\cdot})
$$
where the latter is a simplicial set using the simplicial variable of the
resolution.

We claim that
$$
j\mapsto M_{1/}(F(j), Z_{\cdot})
$$
is a Reedy-fibrant diagram of simplicial sets over $J$. To prove this, note
that the Reedy-fibrant condition is a condition involving strict limits and
colimits, and it is dual to the Reedy-cofibrant condition. Thus
taking $Hom _M= M_{1/}$ of a Reedy-cofibrant diagram $F$, into anything, gives a
Reedy-fibrant diagram as claimed.

Using this claim, and the fact that homotopy limits of simplicial sets may be
calculated using Reedy-fibrant diagrams, we get that
$$
lim ^{\rm str}_{J} M_{1/}(F(j), Z_{\cdot}) \cong
holim _{j\in J}M_{1/}(F(j), Z_{\cdot}).
$$
Here the notation $lim ^{\rm str}_{J}$ means the strict limit taken in the
$1$-category of simplicial sets. On the other hand,
$$
lim ^{\rm str}_{j\in J} M_{1/}(F(j), Z_{\cdot}) =
M_{1/}(colim ^M_JF, Z_{\cdot}).
$$
Putting these all together we get that
$$
L(M)_{1/}(colim ^M_JF, Y) \cong holim _{j\in J}L(M)_{1/}(F(j), Y) .
$$
This exactly says that
$$
colim ^M_JF=colim ^{L(M)}_JF,
$$
meaning in particular that the latter exists.

The case of homotopy limits in $L(M)$ is dual and identical to the above. Thus
we get existence of homotopy limits and colimits in $L(M)$.

Similarly, the cofibrant generation condition implies the second condition in
$(ii)$. Let $\Gg$ be a  small set of objects containing those which occur in a
generating set of cofibrations \cite{Hirschhorn}. Any object of $M$ is weak
equivalent to an object which is obtained as a sequential limit of pushouts by
objects in $\Gg$ (this is the small object argument \cite{Quillen}). Therefore,
$\Gg$ generates $A$ by colimits. Furthermore, the objects of $\Gg$ are
$\beta$-small in $M$ for some $\beta$ (this is in the definition of
``cofibrantly
generated''), and since these objects are cofibrant, $\beta$-smallness in $M$
implies that their images are $\beta$-small in $L(M)$. Thus we obtain
condition (ii).

\subnum{The proof of $(ii) \Rightarrow (iii)$}

Starting with $(ii)$, we can replace our subcategory $C$ which generates $A$ by
colimits, with a subcategory $C$ which strongly generates $A$ by Proposition
\ref{generation}. Therefore we may now assume that the morphism $i: A\rightarrow
\widehat{C}$ is fully faithful.

The next step is to use colimits in $A$ to construct the adjoint functor
$\psi$. Do this as follows. Define the Segal category of arrows
$$
V:= \underline{Hom}(I, \widehat{C}) \times _{\widehat{C}}A
$$
where the first structural morphism in the fiber product is evaluation at $1\in
I$ denoted $ev(1)$.  Thus $V$ is the Segal category of arrows $x\rightarrow y$
with $x\in \widehat{C}$ and $y\in A$. Say that such an arrow is {\em universal}
if it is an initial object in the fiber over $x$ for the evaluation map at $0$
$$
ev(0) : V\rightarrow \widehat{C}.
$$
This condition means that for any object $z\in A$ the morphism
of composition with our arrow,
$$
A_{1/}(y,z) \rightarrow \widehat{C} _{1/}(x,z)
$$
is an equivalence of categories.

Let $U\subset V$ be the full sub-Segal category consisting of universal arrows.

The same definitions can be made with respect to any morphism of
Segal categories $A\rightarrow B$ (in the above we have written
$B=\widehat{C}$).

\begin{lemma}
For any functor $i: A\rightarrow B$, if we construct the Segal category $U$ of
universal arrows from objects of $B$ to objects of $A$, then
the evaluation at $0$ is a fully faithful morphism
$ev(0): U\rightarrow B$.
\end{lemma}
{\em Proof:}
This statement  was shown (in a particular example, but with a technique
which works in general) as a part of the ``second construction'' in the proof of
Lemma 6.4.3 in \cite{SimpsonAspects}. It was this proof which was added in the
revised {\tt v2} of \cite{SimpsonAspects}. We rewrite the proof here.

We may suppose that $A$ and $B$ are fibrant Segal categories. We have
$$
V:= \underline{Hom}(I, B) \times _B A.
$$
This maps by a fibration $ev(0)$ to $B$. Suppose that $u, v\in U\subset V$
are universal arrows; they are maps $u,v:I\rightarrow B$ together with objects
$u_1, v_1$ in $A$ with $u(1)= i(u_1)$ and similarly for $v$. Then
$$
U_{1/}(u,v) = \underline{Hom}(I, B)_{1/}(u,v) \times _{B_{1/}(u(1), v(1))}
A_{1/}(u_1,v_1).
$$
This maps by a fibration to $B_{1/}(u(0), v(0))$. We will show that this latter
map is an equivalence by showing that it is a trivial fibration, namely showing
that it satisfies lifting for any cofibration $E\hookrightarrow E'$.
Thus we suppose that we have a diagram
$$
\begin{array}{ccc}
E & \rightarrow & U_{1/}(u,v)\\
\downarrow && \downarrow \\
E' & \rightarrow & B_{1/}(u(0), v(0)).
\end{array}
$$
The top map amounts to a diagram
$$
\Upsilon (E) \times I \rightarrow B
$$
coupled with a lifting
$$
\Upsilon (E)\times 1 \rightarrow A,
$$
such that over $0$ (resp. $1$) in $\Upsilon (E)$ these restrict to $u$ (resp.
$v$). The bottom map amounts to a diagram
$$
\Upsilon (E') \times 0 \rightarrow B.
$$
We would like to extend the above to a diagram
$$
\Upsilon (E') \times I \rightarrow B
$$
plus lifting
$$
\Upsilon (E')\times 1 \rightarrow A.
$$
Divide the square $\Upsilon (E')\times I$ into two triangles,
$$
\Upsilon (E')\times I = \Upsilon ^2(\ast , E') \cup ^{\Upsilon (E')}
\Upsilon ^2(E', \ast ).
$$
The first triangle is the one containing the edge corresponding to $u$ as its
first edge; the second is the one containing the edge corresponding to $v$
as its
second edge.

We need to define a morphism from these triangles into $B$ (plus a lifting on
the second edge of the first triangle, into $A$). These are already defined
over the subobjects where one puts $E$ instead of $E'$. Furthermore, the first
edge of the second triangle is already defined.

We treat first the second triangle: the inclusion corresponding to the $01$ and
$12$ edges
$$
[\Upsilon (E') \cup ^{\ast} \Upsilon (\ast )]\cup ^{\ldots } \Upsilon ^2(E,
\ast )
$$
$$
\hookrightarrow \Upsilon ^2(E', \ast )
$$
is  a trivial cofibration (cf \cite{SimpsonLimits}) so the extension in
question exists. In particular we obtain an extension along the diagonal.

We now turn to the first triangle $\Upsilon ^2(\ast , E')$. We have an
extension which is already specified along the diagonal (i.e. the $02$ edge)
and we have the specification $u$ which is a universal morphism, along the
$01$ edge. We claim that we can choose the required extension plus lifting
along the $12$ edge into $A$. For this, note (by going back from the
notations $\Upsilon$ to the usual notation) that we are looking at the following
problem. We are given
$$
E' \rightarrow B_{1/}(u(0), v(1)),
$$
and
$$
E \rightarrow \{ u\} \times _{B_{1/}(u(0), u(1))}
B_{2/}(u(0), u(1), v(1))\times _{B_{1/}(u(1),
v(1))} A_{1/}(u_1, v_1),
$$
such that the restriction of
this map to the $02$ edge is the same as the restriction of the first map to
$E$.

We look for an extension
of the above to a map
$$
E \rightarrow \{ u\} \times _{B_{1/}(u(0), u(1))}
B_{2/}(u(0), u(1), v(1))\times _{B_{1/}(u(1),
v(1))} A_{1/}(u_1, v_1).
$$
Recall that
$$
\{ u\} \times _{B_{1/}(u(0), u(1))}
B_{2/}(u(0), u(1), v(1))
\cong B_{1/}(u(1), v(1)).
$$
In particular,
$$
\{ u\} \times _{B_{1/}(u(0), u(1))}
B_{2/}(u(0), u(1), v(1))\times _{B_{1/}(u(1),
v(1))} A_{1/}(u_1, v_1)
\cong A_{1/}(u_1, v_1).
$$
The condition that $u$ is a universal map means that the
$02$-restriction morphism
$$
\{ u\} \times _{B_{1/}(u(0), u(1))}
B_{2/}(u(0), u(1), v(1))\times _{B_{1/}(u(1),
v(1))} A_{1/}(u_1, v_1)
$$
$$
\rightarrow
B_{1/}(u(0), v(1))
$$
is an equivalence. (In other words the ``composition with $u$'' from
$A_{1/}(u_1, v_1)$ to $B_{1/}(u(0), v(1))$ is an equivalence.)
This restriction
map is also fibrant, so it is a trivial fibration and satisfies lifting for
all cofibrations. The lifting condition is exactly the condition that we need
to show.

This completes treatment of the first triangle, and finishes the proof
of the lifting property which shows that the map
$$
U_{1/}(u,v) \rightarrow  B_{1/}(u(0), v(0))
$$
is an equivalence.
\eop

\begin{lemma}
Suppose in the situation of the previous lemma that $B=\widehat{C}$.
If the functor $i : A\rightarrow \widehat{C}$ comes from a morphism
$a:C\rightarrow A$, and if $A$ admits arbitrary small
homotopy colimits, then the morphism  $U\rightarrow \widehat{C}$ is essentially
surjective, so it is an equivalence of Segal categories.
\end{lemma}
{\em Proof:} It suffices to show that if $G\in \widehat{C}$ is any object, then
there exists a universal morphism
$f: G\rightarrow iX$ to an object $X\in ob(A)$. To construct $f$, we
first note that any such $G$ can be expressed as a colimit of objects of $C$
(Lemma \ref{express}):
there is a small category $J$ and a morphism $F: J\rightarrow C$ such that
$$
G= colim ^{\widehat{C}}_J(i \circ a \circ F).
$$
Now we set
$$
X:= colim ^A_J(a\circ F).
$$
We have a morphism of $J$-diagrams in $A$
$$
a\circ F\rightarrow c_J(X),
$$
which gives a morphism of $J$-diagrams in $\widehat{C}$
$$
i\circ a\circ F\rightarrow c_J(iX).
$$
In turn this can be factored through an essentially unique morphism
$f: G\rightarrow iX$ because $G$ is the colimit of $i\circ a\circ F$.
We claim that $f$ is universal. To see this, suppose $Y\in A$.
In what follows we pretend that weak compositions are actually compositions
(this avoids tedious references to things like $A_{2/}$).
We get the following diagram (well-defined and commuting, up to homotopy):
$$
\begin{array}{ccc}
A_{1/}(X, Y) & \rightarrow &
\underline{Hom}(J, A)_{1/}(a\circ F ,c_J(Y))\\
\downarrow && \downarrow \\
\widehat{C}_{1/}(G, iY) & \rightarrow &
\underline{Hom}(J, \widehat{C})_{1/}(i\circ a\circ F ,c_J(iY)).
\end{array}
$$
The top arrow is an equivalence because $X$ is the colimit of $a\circ F$ in $A$.
The bottom arrow is an equivalence because $G$ is the colimit
of $i\circ a \circ F$ in $\widehat{C}$. The right vertical arrow is an
equivalence because of the hypothesis that $i$ is fully faithful (it is the
morphism of functoriality of $i$). Therefore the
left vertical arrow is an equivalence, which exactly says that $f:
G\rightarrow iX$ is a universal arrow.
\eop

Go back to our previous situation (which is the situation of the second half of
the lemma). We have an equivalence $U\rightarrow \widehat{C}$ and
the evaluation at $1$ provides a morphism $U\rightarrow A$; this gives an
essentially well-defined morphism $\psi : \widehat{C}\rightarrow A$.

There is a tautological morphism
$$
U\times I \rightarrow \widehat{C}.
$$
This corresponds to a natural transformation of functors, which (when
we compose with the inverse of the equivalence $U\cong \widehat{C}$) gives a
natural transformation of functors $\widehat{C} \rightarrow \widehat{C}$,
$$
\eta _X: X\rightarrow i\psi (X).
$$
On each object $X$, $\eta _X$ is a universal map.

Now for $Y$ an object of $A$, we look at the morphism induced by $\eta _X$,
$$
A_{1/}(\psi (X), Y) \rightarrow
\widehat{C}_{1/}(X, i(Y)).
$$
The universality condition for $\eta _X$ implies that this map is an equivalence
(of  simplicial sets).  This means that $(\psi , \eta )$ is an adjoint functor
to $i$. This proves condition (iii).

\subnum{Proof of the additional statement in $(iii)$}

In the last paragraph of Theorem \ref{main} was the additional statement that
in (iii), the adjunction morphism going in the other direction is an
equivalence between $\psi \circ i$ and $1_A$. In other words, that $\psi $ is
a retract of the inclusion $i$.

In our discussion of adjunctions, we explained how to obtain the adjunction
morphism in the other direction
$$
\zeta _Y: \psi (i Y) \rightarrow Y.
$$
(Recall, however, that it was left to the reader to show that the morphism
thus constructed was an adjunction.)

We claim that $\zeta _Y$ is an equivalence. To prove this, we will
show that for any $Z\in ob(A)$, the morphism of composition with $\zeta _Y$
$$
A_{1/}(Y,Z) \rightarrow A_{1/}( \psi (i Y), Z)
$$
is an equivalence. Follow this with the two morphisms
$$
A_{1/}( \psi (i Y), Z) \rightarrow
\widehat{C}_{1/}(i\psi (iY), iZ)
\rightarrow
\widehat{C}_{1/}(iY, iZ),
$$
the first of which is functoriality for $i$ and
the second of which comes from the first adjunction morphism $\eta
_{iY}: iY\rightarrow i\psi (iY)$. The composition of these last two morphisms
is an equivalence (that is the condition that $\eta$ be an adjunction). We have
to verify that the composed morphism
$$
A_{1/}(Y,Z)
\rightarrow
\widehat{C}_{1/}(iY, iZ)
$$
is homotopic to the morphism of functoriality for $i$.
For this, note that the diagram
$$
\begin{array}{ccccc}
A_{1/}(Y,Z) &\rightarrow &A_{1/}( \psi (i Y), Z)  &&
\\
\downarrow && \downarrow  && \\
\widehat{C}_{1/}(iY, iZ) &\rightarrow & \widehat{C}_{1/}(i\psi (iY), iZ)
&\rightarrow &
\widehat{C}_{1/}(iY, iZ)
\end{array}
$$
commutes up to homotopy. The bottom row comes from composition with the sequence
$$
iY \stackrel{\eta _{iY}}{\rightarrow}
i\psi (iY) \stackrel{i\zeta _Y}{\rightarrow}
iY.
$$
This composition is homotopic to the identity by Lemma
\ref{adjunctioncompositions}.

Given the above verification,
we get that the composed morphism
$$
A_{1/}(Y,Z) \rightarrow A_{1/}( \psi (i Y), Z)
\rightarrow
\widehat{C}_{1/}(iY, iZ)
$$
is an equivalence (since $i$ is by hypothesis
fully faithful), therefore the first morphism of composition with
$\zeta _Y$ is an equivalence
$$
A_{1/}(Y,Z) \stackrel{\cong}{\rightarrow} A_{1/}( \psi (i Y), Z) .
$$
Since this works for all $Z$, it follows that $\zeta _Y$ is an internal
equivalence in $A$. This proves the claimed statement.

\subnum{The proof of $(iii) \Rightarrow (iv)$}

This part of the proof is mainly concerned with establishing the
Segal-category version of the small generation condition for trivial
cofibrations, which is the main part of the ``cofibrantly generated'' condition
for a closed model category. This benefits from P. Hirschhorn's book
\cite{Hirschhorn} where this
type of condition is widely discussed; and from discussions and
correspondence with A. Hirschowitz about how to put the cardinality arguments
characeristic of Jardine's paper \cite{Jardine}, into a general framework.

Recall from Proposition \ref{generation} that strong generation
implies generation by colimits. Also, the existence of the adjoint $\psi$
implies that $A$ admits all small homotopy colimits (and that $\psi$ preserves
homotopy colimits). In particular, $(iii) \Rightarrow (ii)$ and for the present
step we may use the hypotheses of $(ii)$ and $(iii)$ together. Recall also that
in the previous section of the proof we showed that in the situation of
hypothesis $(iii)$, the adjunction morphism  $\zeta$ is an equivalence between
$\psi \circ i$ and $1_A$.

Let $W\subset \widehat{C}$ be the subcategory of morphisms which go to
equivalences under $\psi$. We obtain morphisms
$$
L(\widehat{C} , W) \stackrel{\psi}{\rightarrow} A,
$$
and
$$
A \stackrel{i}{\rightarrow}L(\widehat{C} , W).
$$
The composition $\psi \circ i$ is already the identity before localization.
On the other hand, the composition $i\circ \psi$ is related to the identity
by a natural transformation which lies in $W$ (we give the argument for this in
the  paragraph which follows); thus on the level of the localization
$i\circ \psi$  is also
homotopic to the identity. This proves that the above morphisms between $A$ and
the localization  $L(\widehat{C} , W)$ are equivalences, and gives the statement
in (iv) about localization.

In the previous paragraph we left open the detail of verifying that for any
$X\in ob(\widehat{C})$, the adjunction morphism $\eta _X : X\rightarrow i\psi
(X)$ is $\psi$-trivial. Thus we have to look at
$$
\psi (\eta _X): \psi (X) \rightarrow \psi (i\psi (X)).
$$
This fits into a sequence
$$
\psi (X) \stackrel{\psi (\eta _X)}{\rightarrow} \psi (i\psi (X))
\stackrel{\zeta _{\psi X}}{\rightarrow}
\psi (X).
$$
By Lemma \ref{adjunctioncompositions}, the
composition is homotopic to the identity of $\psi (X)$. Also the second
morphism is an equivalence as we have shown above. Therefore $\psi (\eta _X)$
is an equivalence, as claimed.

Turn now to our situation
$$
C\rightarrow A \stackrel{i}{\rightarrow} \widehat{C}
$$
where $C$ is a small $1$-category, consists of objects which are small in $A$,
and where $A$ admits small colimits and the morphism $i$ is fully faithful. The
fact that $C$ is small means that the bound for smallness in $A$ of objects of
$C$ can be assumed uniform.  Thus there is a cardinal $\delta$ such that every
object of $C$ is $\delta$-small in $A$.

One has to be careful (cf the ``caution'' after Lemma \ref{transform})
that  the composition $C\rightarrow A\rightarrow \widehat{C}$ is not the Yoneda
morphism for $C$.

By Lemma \ref{objectbyobject}, homotopy colimits in $\widehat{C}$
are calculated object-by-object. Thus we can write
$$
hocolim _{\beta}^S\widehat{C}_{1/}(Z, X_i)\stackrel{\cong}{\rightarrow}
\widehat{C}_{1/} (Z, colim _{\beta }^{\widehat{C}} X_i ) .
$$

We claim that  for any ordinal $\beta$ of size $\geq \delta$,
sequential $\beta$-colimits in $A$ agree with those in $\widehat{C}$. Using the
fact that the objects of $C$ are $\delta$-small in $A$, we get that
for ordinals $\beta \geq \delta$ and sequential $\beta$-diagrams $\{
X_i\}$, we have
$$
hocolim _{\beta}^S\widehat{C}_{1/}(Z, X_i)\stackrel{\cong}{\rightarrow}
\widehat{C}_{1/} (Z, colim _{\beta }^{A} X_i ) .
$$
Comparing with the previous result we get that for any $Z\in C$, the morphism
$$
\widehat{C}_{1/} (Z, colim _{\beta }^{\widehat{C}} X_i )
\rightarrow
\widehat{C}_{1/} (Z, colim _{\beta }^{A} X_i )
$$
is an equivalence. This implies that the morphism
$$
colim _{\beta }^{\widehat{C}} X_i \rightarrow
colim _{\beta }^{A} X_i
$$
is an equivalence in $\widehat{C}$. Thus the two colimits agree.

Now go directly to the case of looking at the functor
$\psi : \widehat{C} \rightarrow A$. Recall that an arrow
$U\rightarrow V$ in $\widehat{C}$ is {\em $\psi$-trivial} if $\psi
U\rightarrow \psi V$ is an equivalence in $A$. An arrow  $F\rightarrow G$ in
$\widehat{C}$ is {\em $\psi$-fibrant} if for every $\psi$-trivial morphism
$U\rightarrow V$, the morphism
$$
Hom (V, F) \rightarrow Hom (V, G) \times _{Hom (U,G)} Hom (U, F)
$$
is an equivalence (where the fiber product is a homotopy fiber product of
simplicial sets). We call this condition the {\em lifting condition}.
A subset $G$ of $\psi$-trivial morphisms is a {\em generating subset} if
a morphism $A\rightarrow B$ which satisfies the above lifting condition
with respect to morphisms  in $G$, is necessarily $\psi$-fibrant.

We would
like to show that the class of $\psi$-trivial morphisms admits a small
generating subset. For this, we adopt the strategy used by Jardine in
\cite{Jardine}. We will take an arbitrary $\psi$-trivial morphism
$U\rightarrow V$ and express it as a sequential colimit of pushouts by
$\psi$-trivial morphisms between smaller objects (until getting back to objects
of a fixed size $\delta$, where we stop and say that we have a generating
subset). To start, we need to be able to talk about the ``size'' of an object.
We say that an object $U\in \widehat{C}$ has {\em size $\leq \beta$} if
$U$ can be expressed (as in Lemma \ref{express}) as a  colimit of objects of
$h_C(C)$, over a $1$-category $J$ with $|J| = \beta$. If this is the case,
then rearrangement of the colimit allows one to express $U$ as a
$\beta$-sequential colimit of objects $U_i$ such that the $U_i$ are of size
$\leq \gamma _i < \beta$.

Fix an ordinal $\delta$ (at least as big as the $\delta$ above, but also at
least uncountable, at least as big as $|C|$, etc.). Let $\Ff$ be the small
set of
$\psi$-trivial morphisms between objects of size $\leq \delta$  in $\widehat{C}$
(technically speaking, this would still be a class but can be replaced by a
small subset containing morphisms equivalent to all those in the class).
We claim that the class of morphisms $\overline{\Ff}$ generated by $\Ff$ in
the sense of lifting, is equal to the full class of $\psi$-trivial morphisms.
(This will complete the proof of $(iii)\Rightarrow (iv)$.)

To prove this, we proceed by transcendental induction: suppose it isn't the
case, and let $\beta$ be the smallest ordinal such that there exists a morphism
$U\rightarrow V$ between objects of size $\leq \beta$, which is $\psi$-trivial
but not in $\overline{\Ff}$. We will show that the morphism may be
expressed as a
sequential colimit of pushouts by smaller $\psi$-trivial morphisms. By the
induction hypothesis and since they are smaller, these $\psi$-trivial morphisms
are in $\overline{\Ff}$; but then this implies that our original morphism is in
$\overline{\Ff}$, a contradiction showing the claim.

Note that $\beta$ is the first ordinal of its cardinality; also $\beta >
\delta$. Express $U$ and $V$ as sequential colimits
$$
U= colim _{\beta} ^{\widehat{C}}U_i
$$
and
$$
V= colim _{\beta} ^{\widehat{C}}V_i,
$$
with the $U_i$ and $V_i$ of size $\leq |i| < \beta$.
The $U_i$ are $\beta$-small in $\widehat{C}$ which implies that, after possibly
reindexing the second colimit, we can assume that the map $U\rightarrow V$
comes from a collection of maps $U_i\rightarrow V_i$.

By assumption, $\psi U \rightarrow \psi V$ is an equivalence in $A$
The fact that $\psi$ preserves colimits means that
the morphisms
$$
colim _{\beta}^A\psi U_i \rightarrow \psi U
$$
and
$$
colim _{\beta}^A\psi V_i \rightarrow \psi V
$$
are equivalences. The
fact that $\beta$-colimits agree in $A$ 	and $\widehat{C}$ means that
the morphisms
$$
colim _{\beta}^{\widehat{C}}\psi U_i \rightarrow \psi U .
$$
and
$$
colim _{\beta}^{\widehat{C}}\psi V_i \rightarrow \psi V .
$$
are equivalences. A similar argument shows that these colimits are
essentially sequential, so (by restricting our attention to big enough indices
$i$) we may assume that they are sequential.

Furthermore by Lemma \ref{objectbyobject} the above colimits in
$\widehat{C}$ are
calculated object-by-object.
Thus, for every $z\in C$ the morphism
$$
hocolim _{\beta}^S (\psi U_i)(z) \rightarrow
hocolim _{\beta}^S (\psi V_i)(z)
$$
is a weak equivalence of simplicial sets.
This implies that there are subsequences $i_k$ and $j_k$ in $\beta$
(which are again indexed by an ordinal which we denote $\kappa$ even though it
is isomorphic to $\beta$), such that
$$
U_{j_k}(z)\rightarrow (\psi V_{i_k})(z)\cup ^{(\psi U_{i_k})(z)} (\psi
U_{j_k})(z)
$$
are equivalences. This fact about simplicial sets comes from Jardine's argument
\cite{Jardine}. Furthermore, since $\beta$ is big with respect to the
cardinality of $C$, we can assume that these same subsequences work for all
$z\in C$. Therefore the morphism
$$
U_{j_k}\rightarrow (\psi V_{i_k})\cup ^{\psi U_{i_k}} (\psi
U_{j_k})
$$
is an equivalence in $\widehat{C}$.

Now apply Lemma \ref{samecolims} which says
that since the $\widehat{C}$-coproduct
$$
(\psi V_{i_k})\cup ^{\psi U_{i_k}} (\psi U_{j_k})
$$
is in $A$ (because it is equivalent by the previous paragraph to $\psi U_{j_k}$)
then this is also the coproduct in $A$. Now the fact that $\psi$ commutes with
colimits means that
$$
(\psi V_{i_k})\cup ^{\psi U_{i_k}} (\psi
U_{j_k})
\cong \psi (V_{i_k} \cup ^{U_{i_k}} U_{j_k}).
$$
Therefore we finally get that the morphism
$$
\psi U_{j_k} \rightarrow \psi (V_{i_k} \cup ^{U_{i_k}} U_{j_k})
$$
is an equivalence. Thus, the morphism
$$
U_{j_k} \rightarrow V_{i_k} \cup ^{U_{i_k}} U_{j_k}
$$
is $\psi$-trivial.
Defining
$$
U'_k := U_{j_k}, \;\;\;
V'_k:= V_{i_k} \cup ^{U_{i_k}} U_{j_k},
$$
we still have $U= colim _{\kappa} U'_k$ and 
$V= colim _{\kappa} V'_k$,
but now the morphisms $U'_k \rightarrow V'_k$ are $\psi$-trivial.
The $U'_k$ and $V'_k$ have size $\leq |j_k| < \beta$, so these $\psi$-trivial
morphisms are in $\overline{\Ff}$. It follows that the morphism
$U\rightarrow V$ is in $\overline{\Ff}$.

\subnum{The proof of $(iv) \Rightarrow (i)$}

Let $N$ be the Heller model category of simplicial presheaves on $C$
\cite{Heller}. Thus
$L(N)\cong \widehat{C}$ (Lemma \ref{heller}).

Let $M$ be the model category with the same underlying
category as $N$, and the same class of cofibrations, but where a morphism  is
said to be a weak equivalence if its image in $\widehat{C}$ is
$\psi$-trivial. As
generating set of trivial cofibrations, we can  take a generating set for $N$
plus a set of cofibrant representatives for our small generating set given in
the hypothesis of (iv). It is easy to see that a morphism is a fibration in $M$
if and only if it is a fibration in $N$, and if its image in $\widehat{C}$ is a
$\psi$-fibration. This implies that the given generating set indeed generates
the trivial cofibrations (that is, lifting for the given generating set is
equivalent to being a fibration i.e. to lifting for all trivial cofibrations).

We use the criterion of \cite{HirschowitzSimpson} Lemma 2.5 to obtain a closed
model structure for $M$ (the numbers in the present paragraph refer to the
conditions in that lemma). As a historical point, note that this lemma is just
a synopsis of the techniques of Dwyer, Kan and Hirschhorn \cite{DHK}
\cite{Hirschhorn}.
Start by noting that $M=N$
admits small limits and colimits (0). Since $M$ is a category of simplicial
presheaves, any small subset is adapted to the small object argument so
conditions (4) and (5) are automatic. The three for two condition (2) is
automatic in view of the definition of weak equivalence. A morphism which
satisfies lifting for all cofibrations is an equivalence in $N$ already so it is
an equivalence in $M$; this gives (3). The cofibrations are the same as for $N$
so condition (6) comes from that of $N$. Condition (7), that the trivial
cofibrations are stable under coproduct and sequential colimit, comes from the
same property for $\psi$-trivial morphisms in $\widehat{C}$. In effect, a
coproduct or sequential colimit of cofibrations, calculated in $M$, is a
homotopy colimit (cf \cite{DHK} \cite{Hirschhorn} \cite{HirschowitzSimpson}), in
other words it is a colimit in $\widehat{C}$; and the $\psi$-trivial
morphisms in
$\widehat{C}$ are stable under coproduct and  sequential colimit because $\psi$
preserves colimits by hypothesis. Finally, for condition (1) note that
cofibrations are stable under retracts because they are the same as for $N$.
As for weak equivalences,  a morphism is by definition a weak equivalence if and
only if it is an equivalence in $A$, and this condition is equivalent to saying
that it projects to an isomorphism in $ho(A)$. The class of isomorphisms in
$ho(A)$ is closed under retracts, so this implies that the class of weak
equivalences in $M$ is closed under retracts. This gives (1).

Therefore by Lemma 2.5 of \cite{HirschowitzSimpson}, we
obtain a cofibrantly generated closed model structure $M$.

To complete the proof of (i) we just have to show that $L(M) \cong A$.
For this, note that $L(M)$ is obtained from $L(N)$ by inverting the
images of $M$-weak equivalences (since $L(N)$ was obtained from $N=M$ by
inverting a subset of the weak equivalences). We have that $L(N)\cong
\widehat{C}$, and the images of the $M$-weak equivalences are exactly the
morphisms of $\widehat{C}$ whose image by $\psi$ is an equivalence in $A$.
The hypothesis of (iv) says that this localization gives exactly $A$.
This completes the proof of $(iv)\Rightarrow (i)$.

We have now finished the proof of Theorem \ref{main}.
\eop

\end{document}